\documentclass[12pt]{amsart}
\usepackage{amssymb}

\input epsf.tex

\def\myitemmargin{
\leftmargini=30pt 
\leftmarginii=21pt 
\leftmarginiii=19pt 
\leftmarginiv=19pt  
}

\newdimen\xsize
\newdimen\oldbaselineskip
\newdimen\oldlineskiplimit
\def\restorelineskip{\baselineskip=\oldbaselineskip%
\lineskiplimit=\oldlineskiplimit}
\def\putm[#1][#2]#3{
\hbox{\vbox to 0pt{\parindent=0pt%
\vskip#2\xsize\hbox to0pt{\hskip#1\xsize $#3$\hss}\vss}}}%
\def\putt[#1][#2]#3{
\vbox to 0pt{\noindent\hskip#1\xsize\lower#2\xsize%
\vtop{\restorelineskip#3}\vss}}

\DeclareFontFamily{U}{rsf}{\skewchar\font'177}%
\DeclareFontShape{U}{rsf}{m}{n}{<-6>rsfs5<6-8>rsfs7<8->rsfs10}{}%
\DeclareFontShape{U}{rsf}{b}{n}{<-6>rsfs5<6-8>rsfs7<8->rsfs10}{}%
\DeclareMathAlphabet\RSFS{U}{rsf}{m}{n}
\SetMathAlphabet\RSFS{bold}{U}{rsf}{b}{n}
\let\scr=\rfs
\def\mib#1{\boldsymbol{#1}}
\def\sf#1{{\mathsf{#1}}} 

\def\slsf{\slshape \sffamily }
\def\msmall#1{\mathchoice{\hbox{\small$\displaystyle {#1}$}}{#1}{#1}{#1}}

\let\3=\ss
\def\cc{{\mathbb C}}

\def\rr{{\mathbb R}}

\def\pp{{\mathbb P}}
  \def\cp{\cc\pp}
\def\zz{{\mathbb Z}}

\def\st{_{\mathsf{st}}}

\def\codim{\sf{codim}\,}

\def\dim{\sf{dim}\,}

\def\dimr{\dim_\rr}
\def\Def{\sf{Def}}

\def\length{\sf{length}}
\def\sym{\sf{Sym}}

\def\sfh{\sf{H}}

\def\hur{{\mib{H}}}
\def\id{\sf{Id}}
\def\im{\sf{Im}\,}

\def\ker{\sf{Ker}\,}
\def\lim{\mathop{\sf{lim}}}

\def\length{\sf{length}}

\def\pr{\sf{pr}}
\def\pt{\sf{pt}}

\def\Sing{\sf{sing}}

\def\supp{\sf{supp}\,}

\def\mbfu{{\mib{u}}}
\def\mbfv{{\mib{v}}}

\def\mbfx{{\mib{x}}}

\def\eps{\varepsilon}

\let\brel=\buildrel
\def\<{\langle}\let\la=\<
\def\>{\rangle}\let\ra=\>
 \let\bs=\bss

\def\d{\partial}

\def\ddef{\mathrel{{=}\raise0.3pt\hbox{:}}}
\def\deff{\mathrel{\raise0.3pt\hbox{\rm:}{=}}}
\def\inv{^{-1}}
\def\hook{\hookrightarrow}
\def\fraction#1/#2{\mathchoice{{\msmall{ #1\over#2}}}%
{{ #1\over #2 }}{{#1/#2}}{{#1/#2}}}
\def\norm#1{\Vert #1 \Vert}

\def\le{\leqslant}

\def\lrar{\longrightarrow}
\def\llar{\longleftarrow}

\def\scirc{\mathop{\mathchoice{\hbox{\small$\circ$}}{\hbox{\small$\circ$}}%
{{\scriptscriptstyle\circ}}{{\scriptscriptstyle\circ}}}}
\def\longpoints{\leaders\hbox to 0.5em{\hss.\hss}\hfill \hskip0pt}
\def\stateskip{\smallskip}
\def\state#1. {\stateskip\noindent{\bf#1. }} 
\def\statep#1. {\stateskip\noindent{\bf#1 }} 

\def\proof{\state Proof. \.}
\def\Chi{\raise 2pt\hbox{$\chi$}}
\def\eg{\hskip1pt plus1pt{\sl{e.g.\/\ \hskip1pt plus1pt}}}
\def\ie{\hskip1pt plus1pt{\sl i.e.\/\ \hskip1pt plus1pt}}
\def\iff{if and only if }
\def\wrt{with respect to }

\def\sli{{\sl i)\/} \.}         \def\slip{{\sl i)}}
\def\slii{{\sl i$\!$i)\/} \.}   \def\sliip{{\sl i$\!$i)}}
\def\sliii{{\sl i$\!$i$\!$i)\/} \.}    \def\sliiip{{\sl i$\!$i$\!$i)}}
\def\sliv{{\sl i$\!$v)\/} \.}

\def\sing{^\sf{sing}}

\def\barr#1{\mskip1mu\overline{\mskip-1mu{#1}\mskip-1mu}\mskip1mu}
\def\Chi{\raise 2pt\hbox{$\chi$}}

\let\phI=\phi\let\phi=\varphi\let\varphi=\phI
   \let\bfla=\bflambda%
\def\scra{\scr{A}}

\def\scrc{\scr{C}}
\def\scrd{\scr{D}}
\def\scrh{\scr{H}}

\def\scrg{\scr{G}}
\def\scri{\scr{I}}
\def\scrj{\scr{J}}

\def\scrm{\scr{M}}

\def\scrn{{\scr N}}
\def\scro{{\scr O}}
\def\scrp{{\scr P}}

\def\scrt{{\scr T}}
\def\scru{{\scr U}}
\def\scrv{{\scr V}}

\def\scry{{\scr Y}}
\def\scrz{{\scr Z}}

\def\eps{\varepsilon}

\def\bs{\backslash}
\def\ogran{{\hskip0.7pt\vrule height6pt depth2.5pt\hskip0.7pt}}

\def\d{\partial}

\def\1{{1\mkern-5mu{\rom l}}}

\def\ge{\geqslant}
\def\inv{^{-1}}
\let\wh=\widehat
\let\wt=\widetilde
\def\fraction#1/#2{\mathchoice{{\msmall{ #1\over#2}}}%
{{ #1\over #2 }}{{#1/#2}}{{#1/#2}}}

\def\le{\leqslant}

\def\.{\thinspace}

\def\ti#1{{\tilde{#1}}}

\def\qed{\ \ \hfill\hbox to .1pt{}\hfill\hbox to .1pt{}\hfill $\square$\par}

\def\comment#1\endcomment{}

 \belowdisplayskip=\abovedisplayskip

\def\lineeqqno(#1){\hfill\llap{\vbox to 10pt%
{\vss\begin{align} \eqqno(#1)\end{align}\vss}}\vskip1pt}

\def\Matrix{\begin{matrix}}
\def\Endmatrix{\end{matrix}}
\def\Cases{\begin{cases}}
\def\Endcases{\end{cases}}

\advance\headheight 1.2pt

\def\newsection[#1]#2{\section{#2}\label{sec:#1}}
   \def\refsection#1{{\slsf Section \ref{sec:#1}$\,$}\showlabel{\tt #1}}
\def\newsubsection[#1]#2{\subsection{#2}\label{sec:#1}\showlabel{\tt #1}}
   \def\refsubsection#1{{\slsf Subsection \ref{sec:#1}$\,$}\showlabel{\tt #1}}

\def\inDex#1{\index{#1}}

\newtheorem{thm}{Theorem}[section]
   \def\newthm#1{\begin{thm} \showlabel{\tt #1}\label{#1}%
      \index{{\slsf Theorem \ref{#1}},  \ label {\tt #1}|Zpage}} 
   \def\refthm#1{{\slsf Theorem \ref{#1}$\,$\showlabel{#1}}}
\newtheorem{lem}[thm]{Lemma}
   \def\newlemma#1{\begin{lem} \showlabel{\tt #1}\label{#1}%
      \inDex{{\slsf Lemma \ref{#1}},  \ label {\tt #1}|Zpage}} 
   \def\lemma#1{{\slsf Lemma \ref{#1}$\,$\showlabel{#1}}}
\newtheorem{prop}[thm]{Proposition}
   \def\newprop#1{\begin{prop} \showlabel{\tt #1}\label{#1}%
      \inDex{{\slsf Proposition \ref{#1}},  \ label {\tt #1}|Zpage}}
   \def\propo#1{{\slsf Proposition \ref{#1}$\,$\showlabel{#1}}}
\newtheorem{corol}[thm]{Corollary}
   \def\newcorol#1{\begin{corol} \showlabel{\tt #1}\label{#1}%
      \inDex{{\slsf Corollary \ref{#1}},  \ label {\tt #1}|Zpage}}
   \def\refcorol#1{{\slsf Corollary \ref{#1}$\,$\showlabel{#1}}}
\newtheorem{defi}{Definition}[section]
   \def\newdefi#1{\begin{defi} \showlabel{\tt #1}\label{#1}\rm %
      \inDex{{\slsf Definition \ref{#1}},  \ label {\tt #1}|Zpage}}

\def\eqqno(#1){\label{eq#1}\tag{\showlabel{#1}}%
    \index{{\slsf  Eqtn \ref{eq#1}},  \ label {\tt(#1)}|Zpage}}
\def\eqqref(#1){\eqref{eq#1}}

\def\showlabel#1{\relax}

\def\eqqno(#1){\label{eq#1}}%
\numberwithin{equation}{section}

\begin{document}

\myitemmargin
\baselineskip =14.0pt plus 2pt

\title[Local Severi problem]%
{On the local version of the Severi problem}
\author[V.~Shevchishin]{Vsevolod V.~Shevchishin}
\address{Fakult\"at f\"ur Mathematik\\
Ruhr-Universit\"at Bochum\\ 
Universit\"atsstrasse 150\\
44780 Bochum\\
Germany}
\email{sewa@@cplx.ruhr-uni-bochum.de}
\dedicatory{}
\date{This version: June 2002. 3d version: March 2002. 
2nd version: March 2001. 1st version: June 2000}
\begin{abstract}
For a given singularity of a plane curve we consider the locus of nodal
deformations of the singularity with the given number of nodes and describe
possible components of the locus. As applications, we solve the local symplectic
isotopy for nodal curves in a neighborhood of a given pseudoholomorphic curve
without multiple components and prove the uniqueness of the symplectic isotopy
class for nodal pseudoholomorphic curves of low genus in $\cp^2$ and $\cp^2 \#
\barr\cp^2$.

\end{abstract}
\maketitle
\setcounter{tocdepth}{2}

\setcounter{section}{-1}
\pagebreak[1]

\newsection[intro]{Introduction}

In the famous {\slsf Anhang F} of his book {\slsf ``Vorlesungen \"uber
algebraische Geometrie''} \cite{Sev}, F.~Severi offered a proof of the statement
that the locus of irreducible plane curves of degree $d$ having the prescribed
number $\nu $ nodes and no other singularities is connected. However, his
argument, which involved degenerating the curve into $d$ lines, is not correct. 
The problem was attacked by several authors, see review of Fulton \cite{Ful},
and the correct proof was given by Harris \cite{Ha}, following original ideas of
Severi.

In this paper we consider the local version of the Severi problem. Let $C^*$ be
a germ of a holomorphic plane curve at the origin $0\in \cc^2$ such that $C^*$
has an isolated singularity at $0$. In particular, $C^*$ can be reducible but
has no multiple components. In this case there exists a versal family $\{ C_s \}
_{s \in \Def(C^*, 0)}$ of deformations of $C^*$ with a non-singular
finite-dimensional base $\Def(C^*, 0)$. Here, as the curve $C^*$ itself, the
family $\{ C_s \} _{s \in \Def(C^*, 0)}$ and the ambient plane $\cc^2$ are
understood in the sense of germs of analytic spaces.

Fix an integer $\nu$. Denote by $\Def^\circ_\nu(C^*, 0)$ the locus of the curves
in $\Def(C^*, 0)$ having exactly $\nu$ nodes and no other singularities, and by
$\Def_\nu(C^*, 0)$ its closure. One can show that $\Def_\nu(C^*, 0)$ consists of
deformations of $C^*$ whose total virtual number of nodes is at least $\nu$. 
Further, let $\delta = \delta(C^*, 0)$ be the virtual number of nodes of $C^*$
at $0$. It is easy to show that $\Def_0(C^*, 0)=\Def(C^*, 0)$, $\Def_\nu(C^*,
0)$ has pure codimension $\nu$ in $\Def(C^*, 0)$, is empty for $\nu>\delta$,
while non-empty and irreducible for $\nu=\delta$.

The subject of the {\slsf local Severi problem} is description of irreducible
components of $\Def_\nu(C^*, 0)$ in the remaining case $0<\nu<\delta$.
The principle result of the present paper is

\state Main Theorem. {\it For $\nu< \delta$, every irreducible component
of $\Def_\nu(C^*, 0)$ contains a nodal curve with $\delta$ nodes}.

More precisely, we show inductively that every irreducible component of
$\Def_\nu(C^*, 0)$ contains a component of $\Def^\circ_{\nu+1}(C^*, 0)$. In
other words, every nodal curve in $\Def^\circ_\nu(C^*, 0)$ can be degenerated
inside the same component of $\Def_\nu(C^*, 0)$ into a nodal curve with exactly
one additional node.

The meaning of {\slsf Main Theorem\/} is that there are no ``unexpected''
components of%
\break 
$\Def_\nu(C^*, 0)$, different from ``expected'' ones obtained by the following
construction. First one deforms $C^*$ into a nodal curve $C'$ with $\delta$
nodes, this is a generic curve in the family $\Def_\delta (C^*, 0)$, and then
smooths $\delta- \nu$ nodes of $C'$. In particular, {\slsf Main Theorem\/}
implies that there exist not more than $\binom{\delta}{\nu}$ irreducible
components of $\Def_\nu (C^*, 0)$. Of course, this bound is very rough. However,
a precise description of components of $\Def_\nu (C^*, 0)$ requires a
description of the action of the monodromy group of $\Def_\delta (C^*, 0)$ on
the set of nodes of the curve $C'$.


\smallskip
Author's motivation for study of the local Severi problem was applications to
the {\slsf symplectic isotopy problem}. It was pointed out in the paper
\cite{Sh} that {\slsf Main Theorem} would imply  the solution of the {\slsf local
isotopy problem for {\it nodal} pseudoholomorphic curves}, which is a version of
the local Severi problem for pseudoholomorphic curves. This result and its
application to the symplectic isotopy problem are presented in \refsection{sym}.

\medskip\noindent
{\slsf Acknowledgments.} The author is strongly indebted to V.~Kharlamov,
St.~Nemirovski, and St.~Orevkov for numerous valuable remarks and suggestions
which helped to clarify the problem. Ph.~Eyssidieux pointed out an error in the
first attempt to the proof, based on completely other ideas. The idea of the
present approach appeared during discussions with M.~Kazarian and
E.~Kudryavtseva. Many other valuable remarks and suggestions were done by
H.~Flenner, G.-M.~Greuel, Ziv Ran, J.-C.~Sikorav, and E.~Shustin.

\newsubsection[skim]{Scheme of the proof}
The main idea of the proof is to trace the ramification locus of the projections
of deformed curves onto a fixed coordinate axis $Oz \subset \cc^2$, $Oz \cong
\cc$. This leads to another deformation problem. In this case too, there exists
a semi-universal family $\{ ( C_s, f_s ) \}$ of pairs ``curve + projection''
with a non-singular finite-dimensional base $\Def(C^* /\Delta)$. Let $\Def_\nu
(C^* / \Delta)$ be the preimage of $\Def_\nu (C^*, 0)$ \wrt the natural
``forgetful map'' $\Def(C^* / \Delta) \to \Def(C^*, 0)$. For a generic $s \in
\Def_\nu (C^* /\Delta )$, denote by $B_\nu(s)$ the branching divisor of the
projection $f_s: C_s \to \Delta$. Then $B_\nu(s)$ depends holomorphically on a
generic $s$, and the family $\{ B_\nu(s) \}$ can be holomorphically extended to
the whole germ $\Def_\nu (C^* /\Delta )$.

The idea of the proof can be now reformulated in study of the loci $\scry_k
\subset \Def_\nu (C^* / \Delta)$ given by the condition ``the multiplicity of
$0\in \Delta$ in $B_\nu (s)$ is at least $k$''. This means that we study
specializations of projections $f_s: C_s \to \Delta$, proceeding successively
along the special strata of the discriminant locus in the space of branching
divisors $B_\nu(s)$. We show inductively that at each step of the specialization
$k \mapsto k+1$, for generic $s \in \scry_{k+1}$, the curve $C_s$ has simple
branchings and nodes lying apart the vertical axis $Ow$, and there is the
following alternative for the structure of $C_s$ at the axis $Ow$:
\begin{itemize}
\item either $C_s$ is non-singular at $Ow$, and in this case we can
proceed to the next inductive step;
\item or $C_s$ has exactly one singular points at the the axis $Ow$, at which
$C_s$ has two branches, both non-singular.
\end{itemize}
We show that in the latter case one can produce exactly one desired node.
Since the total tangency order can not exceed the degree $d$, the latter case
must occur and the inductive procedure will terminate.

\smallskip
Observe that essentially the same construction was used in the Harris'
proof \cite{Ha} of the global Severi problem. Namely, he studied the
varieties $V_{d,g, k}$ of irreducible nodal curves of genus $g$ and degree
$d$ in $\cp^2$ having tangency of order $k$ with some (not fixed) line
$\ell$ at some point $p$. He showed inductively in $k$, that if $V_{d,g,
k}$ is non-empty, then a generic curve in every irreducible component of
$V_{d,g, k}$ has no other incidences and admits a degeneration either into
a generic curve in $V_{d,g, k+1}$, or into a curve having exactly one extra
node. To adapt his point of view to our case, we must simply rotate $\ell$
into the axis $Ow$ and consider the projection $f: C \to \Delta$ from the
infinity point of $Ow$. Note that his condition on the projection $f: C \to
\Delta$ is stronger than our, he requires also a single ramification point of
$f: C \to \Delta$ over $0\in \Delta$.

\baselineskip =14.0pt plus .5pt
\nobreak


\vskip1cm

\tableofcontents

\newsection[1]{Deformation of isolated singularities of plane curves}

\newsubsection[0.1]{Isolated singularities of plane curves}
We recall the standard definitions of the deformation theory, see also
\cite{Pal-1}, \cite{Pal-2}, or \cite{Tju}.

\newdefi{def1.1} An {\slsf isolated singularity of a plane curve} is a germ 
of a curve $(C, 0)$ in $\cc^2$ at the origin $0$, such that $C$ is non-singular
at any $z \ne 0 \in C$. In particular, this means that there are no multiple
components of $C$ at $0$.

A {\slsf deformation} of such an isolated singularity $(C, 0)$ is given by an
analytic map $\pi_S: \scrc_S \to S$ between germs of analytic sets $(\scrc_S,
0)$ and $(S, s_0)$ such that $\pi_S$ is flat and the fiber $\pi_S \inv(s_0)$ is
the germ $(C, 0)$. The germs $(S, s_0)$ and $(\scrc_S, s_0)$ are the {\slsf
base} and the {\slsf total germ} of the deformation, respectively. Such a
deformation $\pi_S: \scrc_S \to S$ is also called a family of deformations of
$(C, 0)$.

Two deformations $(\scrc_S, 0)$ and $(\scrc'_S, 0)$ of $(C, 0)$ with
the same base $(S, s_0)$ are {\slsf isomorphic}, if there exists a germ
biholomorphism $\phi: (\scrc_S, 0) \to (\scrc'_S, 0)$ compatible with
projections $\pi_S: \scrc_S \to S$ and $\pi'_S: \scrc'_S \to S$,
respectively. The notion of as {\slsf isomorphism} of isolated singularities
of plane curves is defined similarly. 

If $\pi_S: \scrc_S \to S$ is a family of deformations of $(C, 0)$, $(T, t_0)$ a
germ of an analytic set, and $\phi: (T, t_0) \to (S, s_0)$ an analytic map, then
$\scrc_T \deff \phi^* \scrc_S \deff \scrc_S \times_S T$ is also a deformation of
$(C, 0)$ \wrt the natural projection $\pi_T: \scrc_T \to T$. In this case
$\pi_T: \scrc_T \to T$ is called the {\slsf pulled-back family} or a deformation
obtained by the {\slsf base change}, and $\phi: (T, t_0) \to (S, s_0)$ is called
the {\slsf base change map}.
\end{defi} 

By G.~Tjurina \cite{Tju} (see also \cite{Don}, \cite{Pal-1} and
\cite{Pal-2}), there exists a semi-universal family of deformations of any
given isolated singularity of an analytic space. In our case we have

\newprop{prop0.1} Let $(C, 0)$ be an isolated singularity of a plane
curve. Then there exists a family $\pi_S: \scrc_S \to S$ of
deformations of $(C, 0)$ with the following properties:

\begin{itemize}
\item[\sli] Any deformation family $\pi_T: \scrc_T \to T$ of $(C, 0)$ is
isomorphic to the pulled-back family $\phi^*\scrc_S \to T$ for an
appropriate base change map $\phi: (T, t_0) \to (S, s_0)$.

\item[\slii] Any morphism $\phi: (S, s_0) \to (S, s_0)$, such that the pulled-back
family $\phi^*\scrc_S \to S$ is isomorphic to $\pi_S: \scrc_S \to S$, is an 
isomorphism.

\end{itemize}

Furthermore, assume that $C^*$ is the zero divisor of the germ of a holomorphic
function $f(z,w)$ at $0\in \cc^2$. Let $\scrt^1(C^*,0) \deff \scro_{\cc^2, 0}/
\bigl(f, \frac{\d f}{\d z}, \frac{\d f}{\d w} \bigr)$. Then $\supp(\scrt^1(C^*,
0)) = \{ 0\}$ and $S$ is smooth $n$-dimensional with $n = \length( \scrt^1(C^*,
0))$.

Moreover, let the germs $\phi_1(z,w), \ldots, \phi_n(z,w)$ of holomorphic
functions generate the basis of $\scrt^1(C^*,0)$ over $\cc$. Set $\Phi(z,w; s_1,
\ldots,s_n) \deff f(z,w) + \sum_i s_i \phi_i(z,w)$, $S \deff (\cc^n, 0)$, and
let $\scrc_S$ be the germ at $0$ of the zero divisor of $\Phi$, equipped with
the projection $\pi_S: \scrc_S \to S$ given by $(z,w; s_1, \ldots,s_n) \mapsto
(s_1, \ldots,s_n)$. Then map $\pi_S: \scrc_S \to S$ is a deformation family of
$(C^*, 0)$ with the desired properties.
\end{prop}

The properties \sli and \slii are {\slsf completeness} and {\slsf minimality} of
the family $\pi_S: \scrc_S \to S$, respectively. Notice also that even if we
deform $(C^*, 0)$ as an abstract complex space, the whole deformation consists
of plane curves.

\newsubsection[1.1]{Deformation of plane curves with projection}  
Now we give an explicit description of deformation of plane curves. Instead of
germs, we shall work with closed analytic subsets in the bi-disc,
apriori with multiplicities.

The following notations are used. $\Delta^2$ denotes the bi-disc with the 
standard complex structure and complex coordinates $(z,w)$, $\pr_1:\Delta^2 
\to \Delta$ is the projection on the first factor. For a complex (\ie 
holomorphic) manifold $X$, compact and with a piecewise smooth boundary $\d X$,
we denote by $\scrh(X)$ the space of holomorphic function which are continuous
up boundary $\d X$. Similar notation $\scrh(C)$ is used in the case when $C$ is
a nodal complex curve. Further, we denote by $\scrh(C, X)$ the space of
holomorphic maps which are continuous up boundary $\d C$ and have image in the
interior of $X$.

Let us start with some standard facts about holomorphic curves in bi-disc.

\newlemma{lem1.1.1} \sli Let $C$ be a holomorphic curve in $\Delta^2$, 
possibly with multiple components. Assume that the projection $\pr_1: C \to
\Delta$ on the first factor $\Delta$ is proper. Then $C$ is the zero
divisor of the uniquely defined unitary Weierstra\3 polynomial $P_f(z,w)
\deff w^d+ \sum_{i=1}^d f_i(z) w^{d-i}$ whose coefficients $f_1(z), \ldots,
f_d(z)$ are bounded holomorphic functions, $f_1(z), \ldots, f_d(z) \in
\scro(\Delta)$.

\slii Let $S$ be a (Banach) analytic set and $F(z, w; s)$ a holomorphic 
function on $\Delta^2 \times S$, such that for every $s\in S$ the projection
$\pr_1$ from zero divisor $C_s$ of $F(z, w; s)$ onto the first factor $\Delta$
is proper. Then $F(z, w; s)$ can be uniquely decomposed into the product $F(z,
w; s)= G(z, w; s) \cdot P(z, w; s)$ where $G(z, w; s)$ is a holomorphic
invertible function on $\Delta^2 \times S$ and $P(z, w; s)$ is a Weierstra\3
polynomial of the form $P(z,w; s) =w^d+ \sum _{i=1} ^d f_i(z; s) w^{d-i}$ whose
coefficients $f_1(z; s), \ldots, f_d(z; s)$ are bounded holomorphic functions on
$\Delta \times S$.
\end{lem}

The result is classical and follows essentially from the Weierstra\3 theorems,
see \eg \cite{Gr-Ha}, Chapter 0. The coefficients $f_i(z) \in \scro(\Delta)$ of
the Weierstra\3 polynomial are used as natural coordinates on the space of
curves in $\Delta^2$ with the cycle topology, when a curve is considered as a
divisor.

\smallskip

\newdefi{def1.1.1} Denote by $\scrz\textsc{}^d(\Delta^2)$ the space of $f=f(z)
=(f_1(z), \ldots, \allowbreak f_d(z)) \in \bigl(\scrh(\Delta)\bigr)^d$, for
which the zero divisor $C_f$ of the Weierstra\3 polynomial $P_f(z,w) \deff w^d+
\sum_{i=1}^d f_i(z) w^{d-i}$ lies in $\Delta \times \Delta(r)$ for some $r=
r(f)<1$ and has no singularities at the boundary. This is a Banach manifold
parameterizing curves $C$ in $\Delta^2$ for which the projection $\pr_1:C \to
\Delta$ is proper and has degree $d$. The curve corresponding to $f \in
\scrz^d( \Delta^2)$ will be denoted by $C_f$. We shall identify $C_f$ with $f$
and write $C_f \in \scrz^d(\Delta^2)$.
\end{defi}

\smallskip
The space $\scrz^d (\Delta^2)$ is too large to work with. We shall replace it by
a finite dimensional moduli space of pairs ``curve + projection''. The latter is
defined by dividing out holomorphic ``slidings'' along vertical fibers in
$\Delta^2$.

\newdefi{def1.1.2} Let $C^* \in \scrz^d(\Delta^2)$ be a holomorphic curve which
has no multiple components and singularities on the boundary, $P_0(z, w)$ its
Weierstra\3 polynomial, and $F_0(z, w)\in \scrh(\Delta^2)$ a holomorphic
function of the form $F_0(z, w) = G_0(z, w) \cdot P_0(z, w)$ with a
non-vanishing $G_0(z, w)\in \scrh(\Delta^2)$. Define the sheaf $\scrt^1(C^* /
\Delta) \deff \scro(\Delta^2) \!\bigm/ \! \bigl( F_0, \frac{\d}{\d w} F_0
\bigr)$, where $( F_0, \frac{\d}{\d w} F_0 )$ states for the ideal generated by
$ F_0$ and its derivative.

The set of the singular points of $C^*$ and the set of critical points of the
projection $\pr_1: C^* \to \Delta$ is called the {\slsf singular set of\/} $\pr:
C^* \to \Delta$ or and denoted by $\Sing(C^*/ \Delta)$. 

\end{defi}

\newlemma{lem1.1a} \sli The support of the sheaf $\scrt^1(C^* / \Delta)$ is the
set $\Sing(C^*/ \Delta)$.

\slii Let $\phi_1(z,w), \ldots, \phi_n(z,w) \in \scrh(\Delta^2)$ be
functions generating a basis of $\scrt^1(C^* / \Delta)$ over $\cc$. Then
every $F (z,w) \in \scrh(\Delta^2)$ sufficiently close to $F_0$ can be
uniquely represented in the form 
\begin{equation}\eqqno(1.1a)
\textstyle
F (z,w) = G(z, w) \cdot \bigl(P_0(z, w+ g(z,w) ) + 
\sum_{i=1}^n s_i  \phi_i(z,w) \bigr )
\end{equation}
with a holomorphic function $ G(z, w)\in \scrh(\Delta^2)$, constants $(s_1,
\ldots, s_n) \in \cc^n$, and a Weierstra\3 polynomial $g(z,w) = \allowbreak
\sum_{i=0}^d w^i b_i(z)$ of degree $d$ with holomorphic coefficients
$b_i(z) \in \scrh(\Delta)$.
\end{lem}

\proof \sli First, let us observe that the definition of $\scrt^1(C^* /
\Delta)$ is independent of the particular choice of the function defining
$C^*$. In particular, $\scrt^1(C^* / \Delta) = \scro(\Delta^2) \!\bigm/ \!
\bigl( P_0, \frac{\d}{\d w} P_0 \bigr)$. Further, it is clear that
$\scrt^1(C^* / \Delta)$ vanishes outside $C^*$ and at regular points of
$C^*$ which are not critical points of the projection $\pr_1: C^* \to
\Delta$. Since $F_0$ must vanish at least quadratically at every singular
point of $C^*$, $\frac{\d}{\d w} F_0$ must also vanish at singular points
of $C^*$. So $\scrt^1(C^* / \Delta)$ is non-trivial at such points. 
Finally, observe that the vertical vector field $\frac{\d}{\d w}$ is
tangent to $C^*$ at critical points of the projection $\pr_1: C^* \to
\Delta$, and hence $\frac{\d}{\d w} F_0$ vanish at such points. This yields 
the first assertion of the lemma.

\slii This assertion will follow from the implicit function theorem
provided we solve the corresponding linearized problem. Differentiating
\eqqref(1.1a) we obtain the equation
\begin{equation}\eqqno(1.1b)
\textstyle
\dot F(z,w) = \dot G(z,w) \cdot P_0(z, w) 
+ G_0(z, w) \cdot \bigl(\frac{\d}{\d w} P_0(z, w)  \cdot \dot g(z,w) + 
\sum_{i=1}^n \dot s_i  \phi_i(z,w) \bigr )
\end{equation}
where dotted symbols state for tangent vectors to the corresponding spaces.
The latter equation is equivalent to 
\begin{equation}\eqqno(1.1c)
\textstyle
\dot F = \dot G \cdot P_0
+ \frac{\d}{\d w} P_0  \cdot \dot g + 
\sum_{i=1}^n \dot s_i  \phi_i ,
\end{equation}
where dotted objects vary in the same Banach spaces as above. Application
of the Weierstra\3' division theorem shows that it is sufficient to
consider the special case where $\dot F $ is a Weierstra\3 polynomial of
degree $d-1$ of the form $\sum_{i=0}^{d-1} w^i a_i(z)$ with holomorphic
coefficients $a_i(z) \in \scrh(\Delta)$. Another application of the
Weierstra\3' division theorem shows that after replacing the functions
$\phi_i(z, w)$ by its remainders $\ti \phi_i(z, w)$ after the division on
$P_0(z, w)$ we obtain an equivalent problem. Observe also that the
remainders $\ti \phi_i(z, w)$ are also Weierstra\3 polynomials of degree
$d-1$. It follows that in a solution of the new problem
\begin{equation}\eqqno(1.1d)
\textstyle
\dot F = \dot G \cdot P_0
+ \frac{\d}{\d w} P_0  \cdot \dot g + 
\sum_{i=1}^n \dot s_i  \ti \phi_i 
\end{equation}
the function $\dot G$ must be  also a Weierstra\3 polynomial of
degree $d-1$.

Now consider \eqqref(1.1d) as a system of linear equations on the coefficients
of Weierstra\3 polynomial $\dot G$ and $\dot g$ so that $\dot F - \sum_{i=1}^n
\dot s_i \ti \phi_i$ is the inhomogeneous part. Then the matrix of coefficients
of the linear system is the Sylvester matrix of the polynomials $P_0$ and
$\frac{\d}{\d w} P_0$, so that its determinant is the resultant of the
polynomials $P_0$ and $\frac{\d}{\d w} P_0$, \ie the discriminant of $P_0$ \wrt
the variable $w$. Let us denote this discriminant by $D(z)$. Then $D(z) \in
\scrh(\Delta)$ and the zero set of $D(z)$ is exactly the projection of the
support of $\scrt^1(C^* / \Delta)$. Since $D(z)$ is not vanishing identically,
it follows the uniqueness of the solution of \eqqref(1.1d) with given $\dot
F$. By the hypotheses of the lemma, for a given $\dot F$ there exists a unique
collection of parameters $(\dot s_1, \ldots, \dot s_n)$ such that $\dot F -
\sum_{i=1}^n \dot s_i \ti \phi_i$ lies in the ideal generated by $P_0$ and
$\frac{\d}{\d w} P_0$. It follows then the solvability of the linear problem
\eqqref(1.1b).
\qed

\newcorol{cor1.2a} \sli The length $n$ of the sheaf $\scrt^1(C^* / \Delta)$
equals to the total vanishing order of the discriminant of the Weierstra\3
polynomial of $C^*$.

\slii The length of the sheaf $\scrt^1(C^* / \Delta)$ is constant under
small deformations of $C^*$.
\end{corol}

\proof \sli We maintain the notation used in the proof of \lemma{lem1.1a}. Let
us apply the elementary ideals theory to the Sylvester matrix of $P_0$ and
$\frac{\d}{\d w} P_0$. Since every ideal of $\scrh(\Delta)$ containing $D(z)$ is
principle, we can bring the Sylvester matrix in the diagonal form, so that the
product of the diagonal elements is $D(z)$. Now it is clear that the minimal
number of the correction terms $\dot s_i\ti \phi(z)_i$ needed to solve
\eqqref(1.1d) with given $\dot F$ is the sum of total vanishing orders of the
obtained diagonal elements.

The second assertion follows from the first one.
\qed

\newcorol{cor1.2b} 
\sli Every curve $C^* \in \scrz^d(\Delta^2)$ is isomorphic to a curve $C \in
\scrz^d(\Delta^2)$ defined by a polynomial.

\slii The deformation space $\Def(C^*/\Delta)$ has natural algebraic structure.
\end{corol}

\proof \sli By \lemma{lem1.1a}, it is sufficient to approximate the Weierstra\3
polynomial $P_0$ of $C^*$ by a polynomial $P$ lying in the ideal generated by
$P_0$ and $\frac{\d}{\d w} P_0$.

\slii By {\slsf Part \slip}, we may assume that $C^*$ is algebraic, \ie the
Weierstra\3 polynomial $P_0(z,w)$ of $C^*$ is a polynomial in the usual
sense. Let $\phi_1(z, w), \ldots, \phi_n(z,w)$ be polynomials inducing a basis
of $\scrt^1(C / \Delta) = \scro(\Delta^2) /\bigl( P_0, \frac{\d}{\d w} P_0
\bigr)$ and $F(z,w; t)$ a {\slsf polynomial\/} in variables $z$, $w$, and $t=
(t_1,\ldots, t_k)$, such that $F(z,w; 0)$ is a defining polynomial for $C^*$. We
assert that the functions $G(z,w)$ and $g(z,w)$ solving the equation
\eqqref(1.1a) with r.h.s. $F(z,w; t)$ are polynomials in variables $z$ and $w$,
and that the dependence of the parameters $s=(s_1, \ldots, s_n)$ and
coefficients of $G(z,w)$ and $g(z,w)$ on $t= (t_1,\ldots, t_k)$ is
algebraic. The first assertion means that the degree of $G(z,w)$ and $g(z,w)$
\wrt variables $z$ and $w$ is bounded uniformly in $t$. This fact
follows from the linearization of \eqqref(1.1a) given by \eqqref(1.1b). The
second assertion is simply reformulation of the fact that \eqqref(1.1a) is a
system of algebraic equations on coefficients.  The corollary follows. \qed

\newdefi{def1.1a} Let $C^* \in \scrz^d(\Delta^2)$ be a curve defined by a
polynomial $P_0(z,w)$. Fix polynomials $\phi_1(z,w), \ldots, \phi_n(z,w)$
generating a basis of $\scrt^1(C^* / \Delta)$. Define $\Def(C^* / \Delta)$
as the germ of $s=( s_1, \ldots, s_n) \in \cc^n$ at $s=0$, $\scrc =
\scrc(C^* / \Delta)$ as the divisor of $P(z,w; s) \deff P_0(z,w) + \sum_i
s_i \phi_i(z,w)$, $\scrc_s \subset \Delta^2$ as the fiber over $s$ of the
projection $\pi_\Def: \scrc \to \Def(C^* / \Delta)$, and $\pr_1: \scrc \to
\Delta$ as the projection on the $z$-disc.
\end{defi}

It follows from \lemma{lem1.1a} that $\Def(C^* / \Delta) \brel\pi_\Def
\over \llar \scrc \brel \pr_1 \over \lrar \Delta$ is a universal deformation
family of the curve $C^*$ equipped with the proper projection onto
$\Delta$. In particular, for another choice of $\phi_1(z,w), \ldots,
\phi_n(z,w)$ we obtain an isomorphic family. As usually, we identify the
germ $\Def(C^* / \Delta)$ with a small neighborhood of $s=0$ in $\cc^n$
representing it.

\newlemma{lem1.3a} Let $C^* \in \scrz^d(\Delta^2)$ be a curve and $\{p_1, \ldots,
p_l \} = \Sing(C^*/\Delta)$ the set of singular points of $\pr_1: C^* \to \Delta$. 
Denote by $C^*_j$ the germ of $C^*$ at $p_j$. Then there exist a natural isomorphism
$\psi: \prod_j \Def(C^*_j / \Delta) \cong \Def(C^* / \Delta)$ and a natural
imbedding $\prod_j \scrc(C^*_j / \Delta) \hook \scrc(C^* / \Delta)$ compatible with
the isomorphism $\psi$ and the projections $\pi_\Def: \scrc(C^*_j / \Delta) \to
\Def(C^*_j / \Delta)$.
\end{lem}

\proof Let $P(z,w;s) =P_0(z,w) + \sum_i s_i \phi_i(z,w)$ be the polynomial
defining a family realizing $\pi_\Def: \scrc(C^* / \Delta) \to \Def(C^* /
\Delta)$. Choose disjoint neighborhoods $U_j$ of $p_j$ which are small bi-discs
with sides parallel to $\Delta^2$, such that $C \cap U_j$ lie in $\scrz^{d_j} (
U_j)$ for the corresponding degree $d_j$. Counting parameters, we conclude that
the restrictions of deformation family $\scrc(C^* / \Delta)$ to $U_j$ induce the
desired isomorphism.
\qed

\medskip
Now let us describe deformation families $\Def(C^*  / \Delta)$ of lower
dimension.

\newlemma{lem1.3b} Let $C^* \in \scrz^d(\Delta^2)$ be a curve, and let $n$
be the dimension of $\Def(C^* / \Delta)$.

\sli If $n=0$, then $C^*$ consists of $d$ disjoint discs and the projection 
$\pr_1: C^* \to \Delta$ is a trivial $d$-sheeted covering. 

\slii If $n=1$, then $C^*$ consists of $d-1$ disjoint discs and the projection 
$\pr_1: C^* \to \Delta$ is a trivial covering on $d-2$ of the discs, and a
$2$-sheeted covering with one simple branching on the remaining disc.

\sliii If $n=2$, then the following cases are possible.

\begin{enumerate}
\item[(a-b)] $C^*$ consists of $d-2$ disjoint discs, the projection $\pr_1: C^*
\to \Delta$ is a trivial covering on $d-4$ of the discs, and a $2$-sheeted
covering with one simple ramification on each of the remaining $2$ discs. The
ramification points can be projected onto $2$ distinct points on $\Delta$
{\rm(case (a))} or onto a single point {\rm(case (b))}.
\item[(c)] $C^*$ consists of $d-2$ disjoint discs, the projection $\pr_1: C^* \to
\Delta$ is a trivial covering on $d-3$ of the discs, and a $3$-sheeted
covering with $2$ simple branching on the remaining disc.
\item[(d)] $C^*$ consists of $d-2$ disjoint discs, the projection $\pr_1: C^* \to
\Delta$ is a trivial covering on $d-3$ of the discs, and a $3$-sheeted
covering with one vertical inflection point on the remaining disc.
\item[(e)] $C^*$ consists of an annulus and $d-2$ discs, the components are
disjoint, the projection $\pr_1: C^* \to \Delta$ is a trivial covering on the
discs and a $2$-sheeted covering with $2$ simple branching on the annulus.
\item[(f)] $C^*$ consists of $d$ discs, $2$ of them meets transversally at one
point, the remaining $d-2$ are disjoint, the projection $\pr_1: C^* \to
\Delta$ is a trivial covering on every disc.
\end{enumerate}
\end{lem}

The proof of the lemma is straightforward. Let us observe that in the cases (a),
(c), and (e) we have 2 simple ramifications over 2 distinct points $z_1$ and
$z_2$ on the disc $\Delta$, whereas the cases (b), (d), and (f), respectively,
correspond to the case when the points $z_1$ and $z_2$ collapse.

\newdefi{def1.1b} Let $C^* \in \scrz^d(\Delta^2)$ be a curve. Set $n\deff \dim
\Def(C^* / \Delta)$. Define\break$\Def_\nu ^\circ (C^* / \Delta)$ as the locus of
those $s\in \Def(C^* / \Delta)$ for which the curve $C_s$ has exactly $\nu$
nodes and no other singularities and the projection $\pr_1: C_s \to\Delta$ has
$n -2\nu$ simple branchings distinct from the projections of nodes. Let $\Def_\nu
(C^* / \Delta)$ be the closure of $\Def_\nu ^\circ (C^* / \Delta)$ in $\Def (C^*
/ \Delta)$.

\slii For $s \in \Def_\nu ^\circ (C^* / \Delta)$, denote by $B_\nu(s)$ the
branching divisor of the projection $\pr_1: C_s \to \Delta$ and by $\sigma_{\nu,
i} (s)$, $i=1,\ldots, n-2\nu$, the $i$-th symmetric polynomial of points of
$B_\nu(s)$. Set $\scrd_\nu(s)(z) \deff z^{n-2\nu} + \sum_{i=1} ^{n-2\nu} (-1)^i
\sigma_{\nu, i} (s) z^{n- 2\nu -i}$, so that $\scrd_\nu(s)(z)$ is the unitary
polynomial in $z$ with zero divisor $B_\nu(s)$. Denote $\scrd_\nu(s) \deff
\bigl( \sigma_{\nu, 1} (s), \ldots \sigma_{\nu, n-2\nu} (s) \bigr)$.
\end{defi}

\newlemma{lem1.1.2} \sli $\Def_\nu (C^* / \Delta)$ is an algebraic
subset of\/ $\Def (C^* / \Delta)$ of codimension $\nu$.

\slii The functions $\sigma_{\nu, i} (s)$, $i=1,\ldots, n-2\nu$, are
holomorphic on  $\Def_\nu ^\circ (C^* / \Delta)$ and extend holomorphically 
on $\Def_\nu (C^* / \Delta)$. 
\end{lem}

\proof First, let us observe that the functions $\sigma_{0, i}$, $i=1, \ldots,
n$, are well-defined and holomorphic on the neighborhood of $C^*$ in the whole
space $\scrz^d (\Delta^2)$. This follows from the construction of the functions
$\sigma_{0, i}$ which is as follows. Starting from a curve $C \in \scrz^d
(\Delta^2)$ close to $C^*$, we take its Weierstra\3 polynomial $P_C(z,w)= w^d +
\sum _i w^{d-i} a_i(z)$; compute the discriminant $D_C(z)$ of $P_C$ with respect
to the variable $w$, this is a polynomial in coefficients $a_i(z)$ of $P_C$; and
then represent the discriminant in the form $D_C(z) = \scrd_C(z) \cdot h_C(z)$
with a non-vanishing holomorphic function $h_C(z) \in \scrh(\Delta)$ and a
unitary polynomial $\scrd _C(z)$ in the variable $z$. Then $\scrd_C(z)$ is the
desired polynomial defining the branching divisor of the projection $\pr_1: C \to
\Delta$ for a {\slsf non-singular} curve $C$. In particular, $\scrd _0 (s) (z) =
\scrd_{C_s} (z)$.

Now observe that for generic $s \in \Def_\nu (C^* /\Delta)$ the polynomial
$\scrd _0 (s) (z)$ has the following structure: it has $n-2\nu$ simple zeros
$z'_1, \ldots, z'_{n-2\nu}$ and $\nu$ double zeros $z''_1, \ldots, z''_\nu$. We
contend that the set of unitary polynomials $p(z)$ of degree $n$ having this
structure is given by a quasi-affine set $A_\nu ^\circ$ in $\cc^n$. To show this
let us consider the map $\sym_\nu: \cc^n \to \cc^n$ associating to each
$n$-tuple $(z_1,\ldots, z_n)$ its elementary symmetric polynomials $\sigma_1
(z_1, \ldots, z_n), \ldots, \sigma_n(z_1, \ldots, z_n)$. Let $A_\nu$ be the
$\sym_\nu$-image of the set given by equations $z_1=z_2$, $z_3 = z_4, \ldots,$
$z_{2\nu-1} = z_{2\nu}$. Since $\sym_\nu$ is algebraic and proper, $A_\nu$ is
Zariski closed in $\cc^n$, and $A_\nu ^\circ$ is Zariski open subset of
$A_\nu$. The latter follows from the fact that the complement $A_\nu \bs A_\nu
^\circ$ describes further incidences among the zeros $z_1, \ldots, z_n$ of the
polynomial $p(z)$, so that it can be defined as the image of a union of
appropriate linear subspaces on $\cc^n$ \wrt the map $\sym_\nu$.

As we shall show later, for $\nu \le \delta(C^*)$ the families $\Def_\nu
(C^*/\Delta)$ are non-empty and contains $C^*$. By {\slsf Lemmas \ref{lem1.3a}}
and {\slsf\ref{lem1.3b}\/}, $\Def_\nu ^\circ (C^*/ \Delta)$ has codimension $\nu$
and $\scrd_0(\Def_\nu ^\circ (C^*/\Delta) ) \subset A_\nu ^\circ$. By
continuity, we obtain $\scrd_0(\Def_\nu (C^*/\Delta) ) \subset A_\nu $. Thus
$\Def_\nu (C^*/\Delta)\subset \scrd_0 \inv (A_\nu)$. Comparing codimension we
conclude that $\Def_\nu (C^*/\Delta)$ is a union of some irreducible components
of $\scrd_0 \inv (A_\nu)$. \lemma{lem1.3b} shows which components of $\scrd_0
\inv (A_\nu)$ belong to $\Def_\nu (C^*/\Delta)$: exactly those ones which meet
$\Def_\nu ^\circ (C^*/\Delta)$.  This proves the first part of the lemma.

\medskip
The second part of the lemma is straightforward.
\qed

\newlemma{lem1.1.3} Let $C^*\in \scrz^d(\Delta^2)$ be a curve represented
by a Weierstra\3 polynomial $P^*(z,w)$ such that $\Sing(C^*/\Delta) = \{ 0 \}$. 
Let $\delta$ be the nodal number of $C^*$ at $0$ and $b$ the number of boundary
components of $C^*$. Then

\begin{itemize}
\item[\slsf a)] For any $\nu$ with $0\le \nu \le \delta$ the space $\Def_\nu (C^*
/\Delta)$ contains $C^*$;

\item[\slsf  b)] for any $\nu$ with $0\le \nu \le \delta$ and $s \in \Def_\nu ^\circ
 ( C^*/ \Delta)$ the normalization $\wt C_s$ of $C_s$ has the Euler characteristic
\begin{equation}\eqqno(1.1.1) 
\chi(\wt C_s)= b + 2(\nu-\delta);
\end{equation}

\item[\slsf c)] $\Def_\delta (C^*/\Delta)$ is  irreducible at $C^*$. In particular,
any two nodal curves $C_0, C_1 \in \Def_\delta ^\circ (C^*/\Delta)$, sufficiently
close to $C^*$, can be connected by a path $C_t$ in $\Def_\delta ^\circ
(C^*/\Delta)$, also close enough to $C^*$.
\end{itemize}
\end{lem}

\proof {\slsf Part a)}. It follows from the hypothesis of the lemma and the
Riemann-Hurwitz formula that $C^*$ consists of $b$ discs passing through origin
$0\in \Delta^2$. Let $C^*_j$, $j=1,\ldots, b$, be the irreducible components of
$C^*$ and $u^*_j : \Delta \to \Delta^2$ their parameterizations, \ie holomorphic
maps such that $u^*_j (\Delta) = C^*_j$. We may assume that the first component
of each $u^*_j$ is given by the formula $\pr_1 \scirc u^*_j(\zeta_j) = \zeta_j
^{d_j}$ where $d_j$ is the degree of the projection $\pr_1: C^*_j \to \Delta$
and $\zeta_j$ is a complex coordinate on the parameterizing disc $\Delta$. 
Making an appropriate approximation, we may additionally assume that every map
$u^*_j : \Delta \to \Delta^2$ is polynomial, so that each $C^*_j \subset
\Delta^2$ extends to a rational affine algebraic curve $u^*_j(\cc) \subset
\cc^2$.

Let $(z,w)$ denote the complex coordinates on $\cc^2$. Consider small
perturbations $u'_j: \cc \to \cc^2$ of the maps $u^*_j: \cc \to \cc^2$ in which
the $z$-component of each $u^*_j$ varies in the space of unitary polynomials of
degree $d_j$ whereas the $w$-component remains unchanged. Let $C'_j \deff
u'_j(\cc) \subset \cc^2$, $j=1,\ldots, b$, be the corresponding rational curves,
and $C' \deff \cup_{j=1} ^b C'_j$ their union. Then, for a generic choice of
$u'_j$, the obtained curve $C'$ must be also generic enough. In particular, the
only singularities of $C'$ are nodes, the projections $\pr_1: C'_j \to \cc_z$ onto
the $z$-axis have only simple branchings, $d_j-1$ for each $C'_j$, the
branchings are disjoint from each other and from the projection of nodes. The
number of nodes of $C'$ is $\delta$, this is essentially the definition of the
virtual number of nodes $\delta(C^*, 0)$. By \lemma{lem1.1a}, $C'$ can be
transformed into a curve $C_s \in \Def_\delta ^\circ (C^* /\Delta)$,
sufficiently close to $C^*$.  To obtain the assertion for arbitrary $\nu =
0,\ldots, \delta$, one needs to smooth $\delta -\nu$ nodes on $C_s$. The needed
construction is provided by \lemma{lem1.3a}. 

\smallskip\noindent
{\slsf Part b)}. Let $n \deff \dim \Def(C^* /\Delta)$. Then $\dim \Def _\nu (C^*
/\Delta) = n-\nu$. On the other hand, it follows from {\slsf Lemmas \ref{lem1.3a}}
and {\slsf\ref{lem1.3b}\.} that for $s \in \Def _\nu ^\circ (C^* /\Delta)$ the
branchings of the projection $\pr_1: C_s \to \Delta$ and the projections of
nodes of $C_s$ form a local coordinate system on $\Def _\nu ^\circ (C^* /
\Delta)$. Thus $n-\nu = b -\chi(\ti C_s) + \nu$ by the Riemann-Hurwitz formula.
Finally note that by the definition of $\delta = \delta(C^*, 0)$ the $s \in \Def
_\delta ^\circ (C^* /\Delta)$ the curve $C_s$ consists of $b$ discs. Thus
$\chi(\ti C_s) =b$ in this case. The assertion follows.

\smallskip\noindent
{\slsf Part c)}. The irreducibility of $\Def_\delta (C^* /\Delta)$ is equivalent
to connectedness of $\Def_\delta ^\circ (C^* /\Delta)$ in a neighborhood of
$C^*$, which is the second assertion of {\slsf Part c)}. Let $C_0, C_1$ be two
curves in $\Def_\delta ^\circ (C^* /\Delta)$ sufficiently close to $C^*$. As it
is shown in the proof of {\slsf Part a)}, we may assume that both $C_0$ and $C_1$
extend to rational affine algebraic curves in $\cc^2$, $\cup_{j=1}^b u_{0,
j}(\cc)$ and $\cup_{j=1}^b u_{1,j}(\cc)$ respectively, where the maps $u_{i,j}:
\cc \to \cc^2$ are polynomial. Moreover, we may assume that $u_{i,j}(\Delta)$
are the components of $C_i$, and that for every $j=1,\ldots,b$ the maps
$u_{0,j}$ and $u_{1,j}$ are close on $\Delta(2)$, \ie $\norm{ u_{0, j} (\zeta)
-u_{1,j} (\zeta) } \ll 1$ for $|\zeta| <2$.

Consider the polynomial maps $u_{\lambda, j}(\zeta) \deff (1-\lambda) u_{0,
j}(\zeta) + \lambda u_{1,j}(\zeta)$ where $\lambda$ varies in the unit disc
$\Delta$. Set $C_\lambda \deff \Delta^2 \cap \bigl( \cup_{j=1}^b u_{\lambda, j}
(\cc) \bigr)$. Then we obtain a family of curves in $\scrz^d(\Delta)$
sufficiently close to $C^*$. It follows that there exists a family $\{
C_{s(\lambda)} \}$ in $\Def(C^* /\Delta)$ such that every $C_{s(\lambda)}$ is
isomorphic to $C_\lambda$. Moreover, the dependence $s(\lambda)$ is holomorphic
since so is the family $\{ C_\lambda \}$. Hence for all by finitely many
$\lambda \in \Delta$ the curve $C_{s(\lambda)}$ lies in $\Def_\delta ^\circ (C^*
/\Delta)$. The lemma follows. \qed

\newsubsection[1.3]{Proof of Main Theorem} An analogue of {\slsf Main Theorem}
for pairs ``curve + projection'' is

\newthm{thm1.3.1} Let $C^* \in \scrz^d(\Delta^2)$ be a curve whose unique singular
point is the origin $0 \in \Delta^2$. Let $\delta = \delta(C^*, 0)$ be the
corresponding virtual number of nodes.

Then for every $\nu< \delta$, every irreducible component of\/ $\Def_\nu (C^* /
\Delta)$ contains a component of $\Def ^\circ _{\nu+1}(C^* / \Delta)$.
\end{thm}

It is obvious that \refthm{thm1.3.1}\/ implies {\slsf Main Theorem}.

\smallskip
Let us explain the ideas lying behind the proof of \refthm{thm1.3.1}. A trivial
but important observation is that the correspondence $\bigl(s \in \Def_\nu
^\circ (C^*/\Delta) \bigr) \brel F_\hur \over \longmapsto \bigl( \pr_1 : C_s \to
\Delta \bigr)$ defines a holomorphic map $F_\hur$ between $\Def_\nu ^\circ
(C^*/ \Delta)$ and the Hurwitz scheme $\hur_{d, m}$ of simply branched coverings
$f: C \to Oz$ over the axis $Oz$ of degree $d$ with $m$ branchings, $m\deff n
-2\nu$.  More precisely, the image of $C_s$ under $F_\hur$ is the trivial
extension of the ramified covering $\pr_1 : C_s \to \Delta$ to the covering
$f_s: \wt C_s \to Oz$.  Obviously, the map $\scrd_\nu : \Def_\nu ^\circ
(C^*/\Delta) \to \cc^m$ factories in the composition $\pi_\hur \scirc F_\hur$,
where $\pi_\hur: \hur_{d,m} \to \cc^m \bs D_m$ is the the natural map
associating to each covering $f: C \to Oz$ its branching divisor and $D_m
\subset \cc^m$ is the discriminant locus.

The map $\pi_\hur: \hur_{d, m} \to \cc^m \bs D_m$ is a non-ramified covering and
its monodromy, the subject of the Hurwitz problem, is understood well enough so
that one can show the following: In the case $m\ge d$, every branched covering
$f: C \to Oz$ of degree $d$ with $m$ simple branchings can be degenerated in that
way that two simple branchings ``collapse'' yielding a node. More precisely, one
uses the monodromy of $\pi_\hur: \hur_{d, m} \to \cc^m \bs D_m$ to attain to a
covering $f: C \to Oz$ which have the same monodromy at two simple branchings, say
$z_1, z_2 \in Oz$, and then contracts $z_1$ with $z_2$ producing the desired
node. Notice that the whole construction can be realized by moving a single
branching of $f: C \to Oz$, say $z_1$, along an appropriate path $\gamma$ in $Oz$
winding around remaining branchings $z_2,\ldots, z_m$. Remark that a similar
argument is used in \cite{G-H-S}.

To realize this construction in our setting, it is necessary to have enough room
for maneuvering with branch points of the covering. However, it is not so, and
the reason for the failure is that the map $s\in \Def_\nu (C^* / \Delta) \mapsto
\scrd_\nu(s) \in \cc^m$ is, in general, not proper.

As an example, let us consider deformation of the ordinary triple point.
Shifting the components of the singularity $C^*$ we obtain three nodes which,
after smoothing, yield six branch points $z_1,\ldots, z_6$ of the projection
$\pr_1: C_s \to \Delta$. It follows from \lemma{lem1.3a} that if $z_1,\ldots,
z_6$ are pairwisely distinct, then any small movement of $z_1,\ldots, z_6$ can
be realized by an appropriate deformation of $\pr_1: C_s \to \Delta$. This means
that the image of $\Def (C^* / \Delta)$ in $\cc^6$ contains an open set.
However, there exists no deformation of $C^*$ for which five branch points, say
$z_2, \ldots, z_6$, collapse and the sixth point $z_1$ remains distinct. Let us
assume the contrary and denote the collapsed points $z_2= \cdots =z_6$ by
$z^*$. The monodromy at the points $z_1$ and $z^*$ in the symmetric group
$\sym_3$ must be a transposition and a product of $5$ transpositions,
respectively. However, since the total monodromy of $\pr_1: C_s \to \Delta$ is
trivial, we must have the same monodromy at $z_1$ and at $z^*$. Thus we conclude
that $C_s$ has two components, say $C'_s$ and $C''_s$, and the projection
$\pr_1: C_s \to \Delta$ is an isomorphism on $C'_s$ and 2 sheeted covering with
2 simple branchings on $C''_s$. Now, remembering meaning of the points $z_1$ and
$z^*$, we see that $C'_s$ and $C''_s$ must have a single intersection point $p$
with intersection index 2, whose projection on $\Delta$ is the point $z^*$. But
this implies that $C'_s$ must be vertical at $p$ and have ramification over
$z^*$, a contradiction.

\smallskip
Since we have no possibility to collapse the branchings of the projection
$\pr_1: C_s \to \Delta$ in the desired way, we study what type of a collapse can
be reached. This is the next idea of the proof. For this purpose we fix an
irreducible component of $\Def_\nu (C^* /\Delta)$ and consider the loci of the
parameters $s$ in this component for which the divisor $B_\nu (s)$ has
multiplicity at least $k$ in the origin $0\in \Delta$ with $k=1,2,\ldots$

\newdefi{def1.3.1} Define
\begin{align}
\Def_{\nu,k} (C^* / \Delta ) \deff& \,\{ s \in \Def_\nu (C^* / \Delta ) : 
\scrd_\nu (s)(z) \text{ is divisible by }z^k\}
\\
\Def^* (C^* / \Delta ) \deff& \,\{ s \in \Def (C^* / \Delta ) : 
\scrd_0 (s)(z) = z^n \}\\
\noalign{\noindent
Besides, we fix an irreducible component $\scry$ of $\Def_\nu (C^* /\Delta)$
through $C^*$ and set} 
\scry_k \deff& \scry \cap \Def_{\nu,k} (C^* / \Delta ) 
\\
\scry^* \deff& \scry \cap \Def^* (C^* / \Delta ) 
\end{align}
\end{defi}

Thus the loci $\Def^* (C^* / \Delta ) \subset \Def (C^* / \Delta )$ and
$\scry^* \subset \scry$ are given by the condition ``{\slsf $B_0 (s)$ is
supported in $0\in \Delta$}''.

\newlemma{lem1.3.2} \sli $\Def_{\nu,k} (C^* / \Delta ) \subset \Def_\nu (C^* /
\Delta)$ and $\scry_k \subset \scry$ are analytic subsets of codimension at most
$k$. 

\slii For $s\in \Def^* (C^* / \Delta )$, the curve $C_s$ consists of discs
and has a unique singular point lying on the axis $Ow$. Moreover, $\Def^* (C^* /
\Delta ) \subset \Def_\delta (C^* / \Delta )$.

\sliii Assume that the generic curve in the family $\scry$  is irreducible. Then
the codimension of $\scry^*$ in $\scry$ is at least $d+1$ except the following
cases when the codimension is $d$:
\begin{enumerate}
\item $C^*$ has a node at $0\in \Delta^2$;
\item $C^*$ consists of two smooth discs which are vertical at the origin $0\in
\Delta^2$. 
\end{enumerate}

\end{lem}

\state Remark. {\slsf Part} \slii of the lemma states that the locus $\Def^* (C^*
/ \Delta )$ consists of curves which have the same topological type of the
singularity as $C^*$.

\proof \sli By definition, $\Def_{\nu,k} (C^* / \Delta ) \subset \Def_\nu (C^* /
\Delta)$ and $\scry_k \subset \scry$ are given by $k$ holomorphic equations
$\sigma_{\nu,i} (s) =0$ with $i= n-2\nu,\, n-2\nu-1, \ldots, n- 2\nu -k+1$.

\smallskip
\slii By the definition of $\scrd_0$, for every $s\in \Def^* (C^* / \Delta )$,
all the singular points of $C_s$ and all the ramification points of the
projection $\pr_1: C_s \to \Delta$ lie on the axis $Ow$. Hence irreducible
components of $C_s$ are discs. Moreover, since $C_s$ is connected, there must be
a unique intersection point of the irreducible components of $C_s$. Indeed,
otherwise one would have either a ramification or a crossing over some $z' \ne 0
\in\Delta$, which contradicts the condition $\scrd_0(s)(z)=z^n$. Thus $C_s$ has
a unique singular point $p^*_s$. This implies that $\delta(C_s, p^*_s) = \delta
= \delta(C^*, 0)$, which means $\Def^* (C^* / \Delta ) \subset \Def_\delta (C^*
/ \Delta )$. 

\smallskip
\sliii We shall compare the codimensions of $\scry$ and $\scry^*$ in $\Def(C^* /
\Delta)$. By definition, $\scry$ has codimension $\nu$ in $\Def(C^* / \Delta)$.
To compute $\codim(\scry^* \subset \Def(C^* / \Delta))$, we first observe that
$\scry^*$ is contained in an irreducible analytic set $\Def_\delta(C^* /
\Delta)$ which has codimension $\delta$ in $\Def(C^* / \Delta)$. So $\codim(
\scry^* \subset \Def(C^* / \Delta)) = \codim(\scry^* \subset \Def_\delta (C^* /
\Delta))+ \delta$. Furthermore, since $\scry^* \subset \Def^* (C^* / \Delta)$
it is enough to estimate $\codim(\Def^* (C^* / \Delta) \subset \Def_\delta (C^*
/ \Delta))$. 

In order to estimate the latter it is sufficient to construct a complex manifold
$\scrv$ with a holomorphic map $f: \scrv \to \Def_\delta (C^* / \Delta)$ such
that $f(\scrv) \ni C^*$, and then estimate $\codim(f\inv(\Def ^* (C^* / \Delta))
\subset \scrv)$. By the universality of $\Def(C^* / \Delta)$, such a map $f:
\scrv \to \Def_\delta(C^* / \Delta)$ corresponds to certain family of
deformations of $C^*$. Recall that the family $\Def_\delta (C^* / \Delta)$
parameterizes those deformations of $C^*$ which consists of discs, and the
number of these discs---denoted by $b$---is the number of irreducible components
of $C^*$.

The desired family is constructed as follows. Let $C^*_i$ be the irreducible
components of $C^*$ and $u^*_i: \Delta \to \Delta^2$ parameterizations of
$C^*_i$. Since every projection $\pr_1: C^*_i \to \Delta$ has a unique branching
at the origin $0\in \Delta$, $u_i^*$ can be chosen in the form $(\zeta^{d_i},
\phi^*_i (\zeta))$, where $d_i$ is the degree of the projection $\pr_1: C^*_i \to
\Delta$. In particular, $\sum_{i=1}^b d_i =d$, the degree of $\pr_1: C^* \to
\Delta$. By \refcorol{cor1.2b} we may assume that every $\phi^*_i$ is
holomorphic in some larger disc $\Delta(r)$, $r>1$. Consider holomorphic maps
$u_i: \Delta(r) \to \cc^2$ given by 
\[
u_i(\zeta)= \big(\zeta ^{d_i} + p_i(\zeta), \phi^*_i(\zeta) + q_i(\zeta)\big)
\]
where $p_i(\zeta)$ is a polynomial of degree $d_i-1$ with zero free term and
$q_i(\zeta)$ is a polynomial of degree at most $1$. Let $\scrv$ be a
sufficiently small ball in the space of the coefficients of the polynomials
$p_i$ and $q_i$. We write $p_{i,\mbfv} = p_{i,\mbfv}(\zeta)$, $q_{i,\mbfv}=
q_{i, \mbfv}(\zeta)$, and $u_{i,\mbfv} = u_{i,\mbfv}(\zeta)$ for the polynomials
and holomorphic maps corresponding to the parameter $\mbfv \in \scrv$. The curve
$C_\mbfv$ is then $\big( \cup_{i=1}^b u_{i,\mbfv} (\Delta(r)) \big) \cap
\Delta^2$. Since the family $\{C_\mbfv\}$ depends holomorphic on $\mbfv\in
\scrv$ it is given by a holomorphically map $f: \scrv \to \Def(C^*/\Delta)$.
Moreover, $f (\scrv) \subset \Def_\delta (C^* / \Delta)$ by \lemma{lem1.1.3}.
Set $\scrv^* \deff f\inv(\Def ^* (C^* / \Delta))$. Then, as it was already
noted, $\codim (\scrv^* \subset \scrv) \le \codim (\Def ^* (C^* / \Delta)
\subset \Def_\delta (C^* / \Delta))$. This means that $\codim( \scry^* \subset
\scry) = \codim( \scry^* \subset \Def (C^* / \Delta)) - \nu \ge \codim( \Def ^*
(C^* / \Delta) \subset \Def_\delta (C^* / \Delta)) + \delta -\nu \ge \codim
(\scrv^* \subset \scrv) + \delta -\nu$.

Estimating $\codim (\scrv^* \subset \scrv)$ we first note that $p_{i, \mbfv}
(\zeta) = 0$ if $\mbfv \in \scrv^*$ since otherwise the projection $\pr_1:
C_{i,\mbfv} \to \Delta$ would have branching outside $0\in \Delta$. This define
a linear subspace $\scrv'$ in $\scrv$ of codimension $\sum_{i=1} ^b (d_i -1) =
d-b$. Perturbation of the free term of the polynomial $q_i$ corresponds a to
vertical shift of $C_{i,\mbfv}$. This means that we must impose further $b-1$
conditions to insure that $C_{i,\mbfv}$ pass through the same point on the axis
$Ow$. Together we obtain a linear subspace $\scrv''$ of $\scrv$ of codimension
$d-1$ parameterizing linear maps $q_i(\zeta) = a_i\zeta$. Since $\scrv^* \subset
\scrv''$, the codimension of $\scry^*$ in $\scry$ is not less than $d$ and
strictly larger $d$ if $\nu \le \delta -2$.

Assume that at least one component of $C^*$, say $C^*_1$, is singular at $0 \in
\Delta^2$. Then the parameterizing map $u^*_1$ has the form $u^*_1(\zeta) =
(\zeta ^{d_1},\, \alpha_l \zeta^l + \alpha_{l+1} \zeta^{l+1} + \cdots)$ with
$d_1 \ge2$ and $l\ge2$. Then for every perturbation of $\mbfv$ in $\scrv''$ by
means of a non-zero linear term $q_1(\zeta) = a_1\zeta$ the component
$C_{i,\scrv}$ must have a node outside the axis $Ow$. This means that in this
case the set $\scrv^*$ is contained in the subspace of $\scrv''$ given by the
condition $q_1(\zeta) =0$, and hence $\codim(\scry^* \subset \scry) \ge d+1$.

It remains to consider the case when $\nu= \delta-1$ and every component of
$C^*$ at $0\in \Delta^2$ is non-singular. Formula \eqqref(1.1.1) and the
irreducibility of a generic curve in $\scry$ imply that the number $b$ of the
components of $C^*$ must be $1$ of $2$. However, the possibility $b=1$ is
excluded since otherwise $C^*$ must consist of a single non-singular component
$0\in \Delta^2$. Thus $C^*$ consists of two non-singular components, $C^*_1$ and
$C^*_2$. Here we must distinguish the following three special subcases according
to the degrees $d_1$ and $d_2$ of the projections $\pr_1: C^*_1 \to \Delta$ and
$\pr_: C^*_2 \to \Delta$, respectively:
\begin{itemize}
\item[(a)] both $d_1$ and $d_2$ equal $1$, \ie  both $C^*_1$ and $C^*_2$
project isomorphically onto $\Delta$;
\item[(b)] $d_1=1$ and $d_2>1$;
\item[(c)] $d_1>1$ and $d_2>1$.
\end{itemize}
Obviously, the subcases (b) and (c) correspond to the subcases (1) and (2) of
the lemma, respectively. In the subcases a) the degree $d$ of $\pr_1: C^* \to
\Delta$ is $2$, the Weierstra\3 polynomial has the form $P(z,w)= w^2 +a(z)w +
b(z)$, and its discriminant is $a^2(z) - 4b(z)$. This implies that the map
$\scrd_0(s)$ is surjective in a neighborhood of the value $s^* \in \Def(C^*/ 
\Delta)$ corresponding to the curve $C^*$, and $\codim(\scry^* \subset \Def(C^*/
\Delta)) = 2\delta$. Thus $\codim(\scry^* \subset \Def_{\delta-1} (C^*/ \Delta))
= \delta +1 \ge d+1=3$ except the case $\delta=1$ when $C^*$ has a nodal
singularity at $0 \in \Delta^2$. This finishes the proof.
\qed

\medskip
Let us give a proof \refthm{thm1.3.1} for the special cases which appear in
\lemma{lem1.3.2} \sliiip.

\newlemma{lem1.3.3} Assume that the curve $C^*$ has two irreducible components
at $0$ which are non-singular. Let $\delta = \delta(C^*, 0)$ be the
corresponding virtual number of nodes.

Then for every $\nu< \delta$, every irreducible component of\/ $\Def_\nu (C^* /
\Delta)$ contains a component of $\Def ^\circ _{\nu+1}(C^* / \Delta)$.
\end{lem}

\proof Consider first the special case when for every component $C_i$ of $C^*$
the projection $\pr_1: C^*_i \to \Delta$ has degree 1. Then $\pr_1: C^* \to
\Delta$ has degree 2 and the Weierstra\3 polynomial of $C^*$ is of the form
$P(z,w) = w^2 + \phi(z) w + \psi(z)$. Hence the discriminant $D(z)$ of $P$ in
the variable $w$ is $D(z)= \phi^2(z) -4\psi(z)$, which is linear in $\psi(z)$.
This implies that the map $\scrd_0: \Def(C^* / \Delta) \to \cc^n$, $n= 2\delta$,
is a biholomorphism on the image. Via zeroes of the discriminant $D(z)$ we have
a complete control on what happens in $\Def(C^* / \Delta)$. In particular, a
curve $C$ lies in $\Def_\nu ^\circ (C^* / \Delta)$ \iff the discriminant $D(z)$
has exactly $\nu$ double zeroes and $2(\delta -\nu)$ simple ones. Moreover, such
a curve $C$ can be holomorphically degenerated into a curve lying in $\Def
_{\nu+1} ^\circ (C^* / \Delta)$. This implies the assertion of the lemma for the
special case.

\smallskip
In the remaining case the curve $C^*$ has two irreducible components at the
origin $0 \in \Delta^2$, both non-singular, such that at least one of them is
vertical at $0$. Choose local holomorphic coordinate $(\ti z, \ti w)$ at the
origin $0 \in \Delta^2$ such that the corresponding projection $\ti \pr_1: C^*
\to \Delta$ has degree 2 at the origin. Then in a neighborhood of the origin
every sufficiently small deformation $C_s$ of $C^*$ is given by the Weierstra\3
polynomial $\wt P_s(\ti z, \ti w) = \ti w^2 + \ti \phi_s(\ti z) \ti w + \ti
\psi_s(\ti z)$. Let $\wt \scrd(s) (\ti z)$ be the polynomial of the degree $n$
in $\ti z$, $n\deff 2\delta$, whose zero divisor is the zero divisor of the
discriminant $\wt D_s(z)= \ti \phi^2_s(z) - 4\ti \psi_s(z)$. Then the
coefficients of $\wt \scrd(s) (\ti z)$ define a holomorphic map $\wt \scrd:
\Def(C^* / \Delta) \to \cc^n$ which has maximal rank at the base point $s^*$
corresponding to $C^*$. One uses the projection $\ti \pr_1: C^* \to \Delta$ to
produce the desired additional node on a curve $C$ from $\Def ^\circ_\nu (C^* /
\Delta)$.
\qed

\medskip\smallskip
Now we proceed to the proof of \refthm{thm1.3.1}. Assume that the singularity of
the curve $C^*$ is not of the type treated in \lemma{lem1.3.3}, \ie that $C^*$
has at least 3 components or that at least one irreducible component of $C^*$
is singular at the origin $0$.

For every index $k=1,\ldots, d$, we fix a decreasing sequence of irreducible
components $\scry'_k$ of $\scry_k$ at $C^*$ so that $\scry = \scry'_0 \supset
\scry'_1 \supset \scry'_2 \ldots\,$

\newprop{prop1.3.3} There exists index $k^* \in \{ 2,\ldots, d\}$ such that:

\smallskip
\sli For $k=0,\ldots, k^*-1$ a generic curve $C$ of the family $\scry'_k$ has
the following structure:
\vskip-18pt\ 
\begin{itemize}
\item $C$ is non-singular at the axis $Ow$;
\item the projection $\pr_1: C \to \Delta$ has branching order $k$ at the origin
$0\in \Delta$ and only simple branchings outside the origin;
\item $C$ is nodal with exactly $\nu$ nodes.
\end{itemize}

\slii A generic curve $C$ of the family $\scry'_{k^*}$ has the following
structure:
\vskip-16pt\ 
\begin{itemize}
\item outside the axis $Ow$, $C$ is nodal with exactly $\nu$ nodes;
\item on the axis $Ow$, $C$ has a unique singular point which either is a node
or consists of two non-singular vertical branches.
\end{itemize}
\end{prop}

It follows from the second assertion of the proposition and \lemma{lem1.3.3}
that a generic curve $C$ from $\scry'_{k^*}$ lies in the family $\Def_{\nu+1}
(C^* / \Delta)$ and is a non-singular point there. This in turn implies
\refthm{thm1.3.1} and hence {\slsf Maim Theorem}.

\smallskip
\proof We use induction in $k$ showing that if for some $k=0,1 \ldots$ the
structure of component $\scry'_k$ is given by \slip, then the structure of
$\scry'_{k+1}$ is given by either \sli or \sliip. We take the case $k=0$ as the
base since $\scry'_0$ has the property \sli of the proposition. Observe also
that the maximal possible branching order of the projection $\pr_1: C \to
\Delta$ over $0\in \Delta$ is $d-1$. Thus for some $k^* \le d$ we must obtain
the case \sliip, and the induction will terminate.

So now we suppose that for some given $k<d$ the component $\scry'_k$ has the
properties listed in \slip. Let $\chi_k$ (resp. $\chi_{k+1}$) denote the Euler
characteristic of (the normalization of) a generic curve from $\scry'_k$
(resp. $\scry'_{k+1}$). It follows from the assumption that $\chi_k = b + 2(\nu
-\delta)$, whereas there are following two possibilities for $\chi_{k+1}$:
$\chi_{k+1} = \chi_k$ and $\chi_{k+1} > \chi_k$. We consider these two cases
separately.

\smallskip\noindent
{\slsf Case $\chi_{k+1} = \chi_k$}. To every curve $C$ from
$\Def_\nu(C^*/\Delta)$ whose normalization $\wt C$ satisfies $\chi(\wt C)=
\chi_k$ we shall associate the following data: The {\slsf zero divisor $Z_C=
\sum_i m_i \zeta_i$} of the composition $\wt C \to C \brel \pr_1 \over \lrar
\Delta$, denoted by $\pr_1: \wt C \to \Delta$ and considered as a holomorphic
function, and the collection of the {\slsf multiplicities} $(m_i)$ of the zero
divisor $Z_C$, defined up to reordering.

Observe that the ramification points of the map $\pr_1: \wt C \to \Delta$ are
exactly those $\zeta_i \in \wt C$ for which $m_i \ge2$, and the branching index
of $C$ at $0\in \Delta$ is $\sum_i (m_i-1)$. Thus $\sum_i (m_i-1)= k$ (resp.\
$\ge k+1$) for a generic curve $C$ in $\scry'_k$ (resp.\ in $\scry'_{k+1}$).
Moreover, the multiplicities $(m_i)$ are the same for two generic curves in
$\scry'_{k+1}$.

Recall that by \refcorol{cor1.2b} a generic curve $C$ in $\scry'_{k+1}$ can be
extended to a holomorphic curve $C^+$ in a larger bi-disc $\Delta(r) \times
\Delta$ with $r>1$, such that the Euler characteristic of the normalization $\wt
C^+$ of $C^+$ is still $\chi_k$. Let $f: \wt C^+ \to \Delta(r) \times \Delta$ be
the composition $\wt C^+ \to C^+ \hook \Delta(r) \times \Delta$. Then every
family $f_s: \wt C^+ \to \Delta(r) \times \Delta$ of sufficiently small
perturbations of $f$ parameterized by $s\in \Delta$ induces a deformation family
$C_s$ of $C$ defined by $C_s \deff f_s(\wt C^+) \cap \Delta^2$. Observe that
under condition $\sum_j (m_j-1) \ge k+1$ on multiplicities the curves $C_s$
remain in $\scry_{k+1}$ and hence in $\scry'_{k+1}$ by irreducibility reason. On
the other hand, for an appropriate choice of the family $f_s$, $s\in \Delta$,
the multiplicities of the curves $C_s$ will satisfy the condition $\sum_j
(m_j-1) = k+1$ for any $s \ne 0$. Thus $\sum_j (m_j-1) = k+1$ for a generic
curve $C$ from $\scry'_{k+1}$. In a similar way one shows that a generic curve
$C$ from $\scry'_{k+1}$ must have the properties \sli of the proposition.

\smallskip\noindent
{\slsf Case $\chi_{k+1} > \chi_k$.} Let $C^\dag$ be a generic curve from
$\scry'_{k+1}$ and $p^\dag_1, \ldots, p^\dag_l \in C^\dag$ the singular points
of $C^\dag$ and the ramification points of the projection $\pr_1: C^\dag \to
\Delta$. Choose sufficiently small bi-discs $\Delta^2_j$ centered at $p_j$ such
that for the curves $C^\dag_j \deff \Delta^2_j \cap C^\dag$ the projections
$\pr_1: C^\dag_j \to \Delta_j$ on the $z$-component are proper. In particular,
$\Delta^2_j$ are mutually disjoint. Then by \lemma{lem1.3a} we obtain a natural
decomposition 
\begin{equation}\eqqno(Def.C.dag)
\textstyle
\Def(C^\dag/\Delta) =\prod_j \Def( C^\dag_j /\Delta_j).  
\end{equation}
More precisely, this should be understand as a natural isomorphism of the germs
(and hence of small neighborhoods) of the spaces at the points corresponding to
the curve $C^\dag$.

It follows from \refcorol{cor1.2a} that we can consider the space $\Def( C^\dag
/ \Delta)$ as an open subset of $\Def(C^*/\Delta)$. In particular, the loci
$\scry'_k \cap \Def( C^\dag / \Delta)$ and $\scry' _{k+1} \cap \Def(C^\dag /
\Delta)$ describe the behavior of curves of the families $\scry'_k$ and $\scry'
_{k+1}$ near $C^\dag$. We contend that the decomposition \eqqref(Def.C.dag) is
compatible with the families $\scry'_k$ and $\scry' _{k+1}$.

To show this let us take a generic curve $C$ in $ \Def_{\nu,k} (C^\dag/\Delta)$
sufficiently close to $C^\dag$. Considering the pieces $C_j \deff \Delta^2_j
\cap C$ of $C$, we can ``decompose'' the numerical invariants characterizing
$\Def_{\nu,k} (C^* / \Delta )$. Namely, we obtain the following decompositions:
\begin{itemize}
\item[(a)] $\nu = \sum_j \nu_j$ where $\nu_j$ is the number of
nodes of $C_j$;
\item[(b)] $k = \sum_j k_j$ where $k_j$ is the branching degree of the projection
$\pr_1: C_j \to \Delta_j$ over $z=0$ if the point $p_j$ lies on the $Ow$-axis
and $k_j=0$ otherwise.
\end{itemize}
Using this we obtain further four natural decompositions:
\begin{itemize}
\item[(c)] $\Def_\nu(C^\dag/\Delta) = \bigcup\limits_{\scriptscriptstyle 
\sum_j \nu_j =\nu} \prod_j \Def_{\nu_j} (C^\dag_j /\Delta_j)$;
\smallskip
\item[(d)] $\scrd_\nu(s)(z) = \prod_j \scrd_{\nu_j}(s_j)(z)$ \ for \ 
$s=(s_j)\in \prod_j \Def_{\nu_j} (C^\dag_j /\Delta_j)$;
\medskip
\item[(f)] $\Def_{\nu,k} (C^\dag/\Delta) = \mathop{\bigcup'}\limits_{
{\sum_j\nu_j =\nu \atop \sum_j k_j =k} }
\prod_j \Def_{\nu_j, k_j} (C^\dag_j /\Delta_j)$,
\medskip
\item[(g)] $\Def_{\nu,k+1} (C^\dag/\Delta) = \mathop{\bigcup'}\limits_{
{\sum_j\nu_j =\nu \atop \sum_j k_j =k+1} }
\prod_j \Def_{\nu_j, k_j} (C^\dag_j /\Delta_j)$,
\end{itemize}
where the union in (f) and (g) is made only over those decompositions $k=\sum_j
k_j$ or $k+1=\sum_j k_j$, respectively, which can appear in (b), \ie for which
the component $k_j$ is zero if $p^\dag_j$ does not lies on the axis $Ow$.

Decomposition (c) follows from the definition of the families $\Def_\nu$ and
\lemma{lem1.3a}, decomposition (d) from \lemma{lem1.1.2}, whereas decompositions
(f) and (g) from (d) and also \lemma{lem1.1.2}.

Since $\scry'_k$ and $\scry'_{k+1}$ are irreducible, there exist uniquely
defined decompositions $\nu = \sum_j \nu^\dag _j$, $k= \sum_j k^\dag_j$ and $k+1
= \sum_j k^\ddag_j$, such that $\scry'_k \cap \Def( C^\dag / \Delta)$ lies in
$\prod_j \Def_{\nu^\dag_j, k^\dag_j} (C^\dag_j /\Delta_j)$ and $\scry'_{k+1}
\cap \Def( C^\dag / \Delta)$ lies in $\prod_j \Def_{\nu^\dag_j, k^\ddag_j}
(C^\dag_j /\Delta_j)$. Moreover, there exists the unique index $j_0$, say
$j_0=1$, such that $k^\ddag_j= k^\dag_j$ for $j\ne j_0=1$ and $k^\ddag_1 =
k^\dag_1 +1$.  Observe that the corresponding point $p^\dag _1$ lies on the axis
$Ow$.

The condition of genericity of $C^\dag$ in $\scry'_{k+1}$ implies that $\scry'
_{k+1}$ is non-singular at $C^\dag$. Thus every family $\Def_{\nu^\dag_j,
k^\ddag _j} (C^\dag_j /\Delta_j)$ is non-singular and generic at $C^\dag_j$. This
means that $\Def_{\nu^\dag_j, k^\dag _j} (C^\dag_j /\Delta_j)$ are non-singular
and generic at $C^\dag_j$ for every $j \ne 1$. Consequently, there exists
irreducible components $\scry^\dag_{k^\dag_1}$ of $\Def_{\nu^\dag_1, k^\dag _1}
(C^\dag_1 /\Delta_1)$ and $\scry^\dag_{k^\dag_1+1}$ of $\Def_{\nu^\dag_1, k^\dag
_1+1} (C^\dag_1 /\Delta_1)$ at $C^\dag_1$ such that 
\begin{align}
\eqqno(dec-k)
\scry'_k \cap \Def_{\nu, k} (C^\dag /\Delta) & \textstyle = 
\scry^\dag_{k^\dag_1} \times 
\prod_{j>1} \Def_{\nu^\dag_j, k^\dag_j} (C^\dag_j /\Delta_j)
\\
\eqqno(dec-k+1)
\scry'_{k+1} \cap \Def_{\nu, k} (C^\dag /\Delta) & \textstyle = 
\scry^\dag_{k^\dag_1+1} \times 
\prod_{j>1} \Def_{\nu^\dag_j, k^\dag_j} (C^\dag_j /\Delta_j)
\end{align}
Moreover, every factor in \eqqref(dec-k) represents a family of curves which
satisfies the conditions listed in the part \sli of the hypothesis of the
proposition. 

\smallskip
We contend that the decompositions \eqqref(dec-k) and \eqqref(dec-k+1) are
non-trivial in the sense that the dimension of every factor $\Def_{\nu^\dag_j,
k^\dag _j} (C^\dag_j /\Delta_j)$ is positive and the number of the
factors---which is the number $l$ of the points $p^\dag_j$---is at least 2. The
latter follows from \lemma{lem1.3.2}, \sliiip. Moreover, we have shown that at
least one point $p^\dag_j$ lies not on the axis $Ow$. This implies that the
dimension of every factor $\Def_{\nu^\dag_j, k^\dag_j} (C^\dag_j /\Delta_j)$ is
strictly less that the dimension of $\scry'_k$.

This provides that now we can use the induction in the dimension of the family
$\scry'_k$. This means that since the dimension of $\scry^\dag_{k^\dag_1}$ is
strictly less that the dimension of $\scry'_k$ and since the family
$\scry^\dag_{k^\dag_1}$ has the properties listed in the part \sli of the
proposition, the family $\scry^\dag_{k^\dag_1+1}$ must be either of type \sli or
type \slii of the proposition. In view of properties of the decomposition
\eqqref(dec-k+1), the same dichotomy holds also for $\scry'_{k+1}$. \qed

\medskip
\newsection[sym]{Application to the symplectic isotopy problem} 

\newsubsection[sym1]{The symplectic isotopy problem for nodal surfaces}
We consider a version of the symplectic isotopy problem for surfaces in a
symplectic {\sl 4}-manifold with positive ordinary double points. As an
introduction to the problem we refer to author's paper \cite{Sh}.

\newdefi{def2.1} Let $(X,\omega)$ be a symplectic 4-manifold. A {\slsf
nodal symplectic surface} in $X$ is an immersed surface $\Sigma \subset X$ such
that the restriction $\omega\ogran_\Sigma$ never vanishes and the only
singularities of $\Sigma$ are positive ordinary double points, called {\slsf
nodes}. 

Note that the restriction $\omega\ogran_\Sigma$ induces the orientation on
$\Sigma$. Recall that an ordinary double point of an immersed oriented surface
in a 4-fold is {\slsf positive} if the self intersection number at this point is
$+1$.

Two closed nodal symplectic surfaces $\Sigma_0$ and $\Sigma_1$ in $(X,\omega)$
are {\slsf symplectically isotopic} if they can be connected by an isotopy
$\Sigma_t$ consisting of nodal symplectic surfaces. Such an isotopy $\Sigma_t$
is called a {\slsf symplectic isotopy} between $\Sigma_0$ and $\Sigma_1$.
\end{defi}

Now the {\slsf symplectic isotopy problem} can be formulated as follows:

\smallskip\noindent
{\it Given a symplectic {\sl 4}-manifold $(X,\omega)$ and closed irreducible
nodal symplectic surfaces $\Sigma_0$, $\Sigma_1 \subset X$ lying in the same
integer homology class and having the same genus $g$, does there exists a
symplectic isotopy between $\Sigma_0$ and $\Sigma_1$\sl?}

\smallskip
Note that the genus of a closed irreducible nodal symplectic surface $\Sigma$ in
a symplectic 4-manifold $(X,\omega)$ can be computed by the {\slsf genus formula}
\begin{equation}\eqqno(2.1)
g(\Sigma) = \frac{[\Sigma]^2 -c_1(X,\omega)\cdot[\Sigma]}{2} +1 -\delta(\Sigma),
\end{equation}
where $\delta(\Sigma)$ is the number of nodes on $\Sigma$ and $c_1(X,\omega)$ is
the first Chern class of $(X,\omega)$ (see \eg \cite{Gro}, \cite{McD-Sa-1}, or
\cite{Sh}). Thus in the situation of the symplectic isotopy problem the number
of nodes on $\Sigma_0$ and $\Sigma_1$ is the same.

\smallskip
In the paper \cite{Fi-St} Fintushel and Stern exhibited a class of symplectic
4-folds $(X,\omega)$ with the following property. There exists an infinite
number of symplectic imbeddings $\Sigma_i \hook X$, such that all $\Sigma_i$ are
homologous but pairwise non-isotopic, even smoothly.  So the answer to the
symplectic isotopy problem can be negative in general. On the other hand, the
results of the paper \cite{Sh} give reason to hope that the answer might be
positive for special symplectic 4-folds. Namely, in \cite{Sh} the author
formulated the following

\smallskip
\state Conjecture. {\it Let $\Sigma_0$ and $\Sigma_1$ be closed irreducible
nodal symplectic surfaces in a closed symplectic {\sl 4}-manifold $(X,\omega)$
lying in the same integer homology class and having the same genus $g$. Then a
symplectic isotopy between $\Sigma_0$ and $\Sigma_1$ exists provided
$c_1(X)\cdot [\Sigma_0]>0$.}

\medskip
As it was mentioned {\slsf Introduction}, the solution of the local Severi
problem in the form of {\slsf Main Theorem} implies a solution of the local
symplectic isotopy problem for the case of immersed surfaces with (positive)
nodes. In order to explain this relation, let us make an overview of the method
used for constructing a symplectic isotopy.

\medskip
First, recall that there exists a complete classification of compact symplectic
4-folds $X$ which come in question. 

\newprop{prop2.1}
Let $(X,\omega)$ be a compact symplectic $4$-fold and $\Sigma \subset X$ a
closed symplectic nodal surface with $\la c_1(X), [\Sigma] \ra >0$. Assume that
$\Sigma$ is not an exceptional sphere. Then $X$ is either $\cp^2$, or a ruled
complex surface, or its blow-up.
\end{prop}

For the precise description of the blow-up procedure in symplectic category we 
refer to \cite{McD-3} and \cite{Gi-St}. 

\proof For the case of {\slsf imbedded\/ $\Sigma$}, this proposition is proved in 
\cite{McD-Sa-2}, {\bf Corollary 1.5}. The general case follows from the fact
that every symplectic nodal surface $\Sigma$ in a symplectic $4$-fold can be
``symplectically smoothed'', \ie deformed into an imbedded symplectic surface.
\qed

\medskip
The complete description of possible symplectic structures on such $X$ was given
in \cite{McD-2}, \cite{La-McD}, and \cite{McD-Sa-2}, see also \cite{Li-Liu},
\cite{Liu}.

\newprop{prop2.2} \sli Every symplectic form $\omega$ on $\cp^2$ is 
isotopic to a multiple of the Fubuni-Study form $\omega\st$.

\slii Every symplectic form $\omega$ on a (minimal) ruled complex surface $X$
is compatible with some genuine complex structure $J$.
\end{prop}

The minimality is understood in the sense of ruled complex surfaces so that $X$
is not a blow-up of another ruled complex surface. Thus the $\cp^2$ blown-up
once is minimal in this sense. The compatibility of $J$ and $\omega$ means that
they define a K\"ahler structure on $X$.

\medskip
Now we recall main features of Gromov's theory of pseudoholomorphic curves which
is for the moment the most effective approach to the symplectic isotopy
problem.

\newdefi{def2.2}
An {\slsf almost complex structure} on a manifold $X$ is an endomorphism $J$ of
the tangent bundle $TX$ such that $J^2 = -\id$.  The pair $(X,J)$ is called an
{\slsf almost complex manifold}.

An almost complex structure $J$ on a symplectic manifold $(X,\omega)$ is called
{\slsf $\omega$-tame} if $\omega(v,Jv) \allowbreak >0$ for any non-zero tangent
vector $v$. The set of $\omega$-tame almost complex structures on $X$ is denoted
by $\scrj_\omega$.
\end{defi}

\newdefi{def2.3} A {\slsf parameterized $J$-holomorphic curve} in an almost 
complex manifold $(X,J)$ is given by a Riemann surface $S$ with a
complex structure $J_S$ on $S$ and a (non-constant) $C^1$-map $u: S \to X$ 
satisfying the {\slsf Cauchy-Riemann equation}
\begin{equation}
 du + J \scirc du \scirc J_S =0.
\eqqno(2.2)
\end{equation}
In this case we call $u$ a $(J_S, J)$-{\slsf holomorphic map}, or simply
$J$-{\slsf holomorphic map}. Here we use the fact that if $u$ is not constant,
then the structure $J_S$ is unique. In particular, such a map $u$ equips $S$ 
with a complex structure $J_S$.

A {\slsf non-parameterized $J$-holomorphic curve} is the image $C = u(S)$ of a
non-constant $J$-holomorphic map $u: S \to X$. Since the map $u$ equips $S$ with
a complex structure $J_S$ we obtain a Riemann surface $(S, J_S)$ which can be
seen as the {\slsf normalization} $\wt C$ of $C=u(S)$ provided $C$ is non-multiple.

The structure of $J$-holomorphic maps and curves is very similar to that of
usual holomorphic objects, for details see \eg \cite{Mi-Wh}, \cite{Sk-1},
\cite{Sk-2}, and \cite{Sh}. In particular, the notions of an {\slsf irreducible
component} and the {\slsf multiplicity} of a component have the usual meaning.

The notion {\slsf $J$-holomorphic curve} or simply even {\slsf $J$-curve} means
either parameterized or non-parameterized curve. We say about {\slsf
pseudoholomorphic} maps and curves if the structure $J$ is clear from the
context or not specified.

We always assume that the parameterizing surface $S$ is compact but not
necessary closed, so that the  boundary $\d S$ of $S$ can be non-empty. In this
case we assume that $\d S$ consists of finitely many smooth circles and that
both the structure $J_S$ and the parameterizing map $u$ are $C^1$-smooth up to
boundary $\d S$. The {\slsf boundary $\d C$} of a pseudoholomorphic curve $C$
parameterized by $u: S \to X$ is the set $u(\d S)$. We say that a curve $C$ is
{\slsf non-singular} at the boundary $\d C$ if $u$ is an imbedding near $\d S$.
\end{defi}

Applying Gromov's theory to the symplectic isotopy problem, one uses the
following argumentation. It is well-known that the set $\scrj_\omega$ of tame
almost complex structures in a symplectic manifold $(X,\omega)$ is non-empty and
contractible (see \eg \cite{Gro}, \cite{McD-Sa-1}). In particular, any two
$\omega$-tame almost complex structures $J_0$ and $J_1$ can be connected by a
homotopy (path) $J_t$, $t\in [0,1]$, inside $\scrj_\omega$.  Furthermore, every
immersed surface $\Sigma$ in a symplectic 4-fold $(X,\omega)$ with ordinary
double points is $J$-holomorphic curve with respect to some $\omega$-tame
structure $J$ \iff $\Sigma$ is a nodal $\omega$-symplectic surface.

\smallskip
Now let $(X, J_1)$ be a (compact) ruled complex surface with a K\"ahler form
$\omega$ and $\Sigma$ a nodal $\omega$-symplectic closed surface in $X$.  Find
an $\omega$-tame almost complex structure $J_0$ making $\Sigma$ a
$J_0$-holomorphic curve. Find a path $h: [0,1] \to \scrj_\omega$ such that $h(0)
= J_0$ and $h(1) = J_1$, so that $J_t \deff h(t)$ is a homotopy between $J_0$
and $J_1$. Fix points $\mbfx= (x_1, \ldots, x_k)$ on $X$ and consider the spaces 
\begin{align}
\eqqno(2.3)
\scrm_{h, \mbfx} &\deff \left\{ (C,t) : \;
\vcenter{
\vbox{ \hsize = 0.51\hsize \parindent=0pt
\small
$t\in [0,1]$, $C$ is a non-multiple irreducible $h(t)$-ho\-lomorphic curve of
geometric genus $g$ in the homo\-logy class $[\Sigma]$ passing through $x_1,
\ldots, x_k$ }} \;
\right\},
\\
\eqqno(2.4)
\scrm_{h, \mbfx} ^\circ &\deff \big\{ (C,t) \in  \scrm_{h, \mbfx} :\; 
\text{ $C$ is nodal}\, \big\}
\end{align}
together with the projection $\pr_{h, \mbfx}: \scrm_{h, \mbfx} \to [0,1]$. 

The reason for introducing the points $x_1,\ldots, x_k$ will be explained later.
For a while, we may assume that $k=0$ and there is no constrain on curves to
pass through given points.

It is known that for a {\slsf generic} path $h: [0,1] \to \scrj_\omega$ the
space $\scrm_{h, \mbfx}$ has a natural structure of a smooth manifold of the
expected dimension
\[
\dimr \scrm_{h, \mbfx} = 1 + 2(c_1(X)\cdot [\Sigma] + g-1 -k)
\]
such that the projection $\pr_{h, \mbfx}$ is smooth, and $\scrm_{h, \mbfx}
^\circ$ is open in $\scrm_{h, \mbfx}$. Let us denote by $\pr_{h, \mbfx} ^\circ$ 
the restriction of $\pr_{h, \mbfx}$ onto $\scrm_{h, \mbfx}
^\circ$.

A crucial observation is that a section $s(t) = (C_t,t)$ of the projection
$\pr_{h, \mbfx} ^\circ$ with $C_0 =\Sigma$, if exists, would give a symplectic
isotopy between $\Sigma$ and a holomorphic curve $C_1$. Furthermore, since the
moduli space of nodal $J_1$-holomorphic (and hence algebraic) curves of the
given geometric genus $g$ and homology class $[\Sigma]$ in $X$ is
quasi-projective, it has finitely many components. This would reduce the
symplectic isotopy problem to the {\slsf Severi problem} of $(X, J_1)$: the
description of components of the space $\scrm_{J_1}$ of nodal irreducible curves
in $(X, J_1)$ of given homology class and genus. The case of primary interest
for the symplectic isotopy problem is the one with $c_1(X) \cdot [C] >0$. There
is a certain progress in this direction after Harris' paper, see \cite{Ran} and
\cite{G-L-Sh}. However, the answer to the Severi problem in the case
$c_1(X) \cdot [C]>0$ in general is still unknown.

\smallskip
Constructing of a section $s(t)$ of the projection $\pr_{h, \mbfx}$ one
challenges two principal difficulties. The first one is that the projection
$\pr_{h, \mbfx}: \scrm_{h, \mbfx} \to [0,1]$ considered as a real function can
have local maxima. However, as it was shown in {\slsf Section 4} of \cite{Sh},
this difficulty does not occur if $c_1(X)\cdot [\Sigma] > 0$. More precisely, it
is proved that 
\begin{itemize}
\item[\slsf S1)\.] the complement of $\scrm_{h, \mbfx} ^\circ$ in $\scrm_{h, 
\mbfx}$ has Hausdorff codimension $\ge2$:
\item[\slsf S2)\.] if the number $k$ of fixed points $\mbfx$ is strictly less 
than $c_1(X)\cdot [\Sigma]$, then (for a generic $h$) every critical point of
the projection $\pr_{h, \mbfx}: \scrm_{h, \mbfx} \to [0,1]$ is {\slsf saddle}.
\end{itemize}
This insures that a section of $\pr_{h, \mbfx}$ over $[0,t_0]$ can be continued
to a bigger interval $[0, t_1)$, $t_1>t_0$.

The second difficulty comes from the fact that the space $\scrm_{h, \mbfx}$ is
not compact and the projection $\pr_{h, \mbfx}: \scrm_{h, \mbfx} \to [0,1]$ is
not proper. Gromov's compactness theorem provides that there exists a nice
compactification $\barr\scrm_{h, \mbfx}$ of $\scrm_{h, \mbfx}$ such that 
\begin{itemize}
\item $\barr\scrm_{h, \mbfx}$ is a compact Hausdorff topological space;
\item it has a natural stratification whose strata are smooth for a generic $h$;
\item $\pr_{h, \mbfx}: \scrm_{h, \mbfx} \to [0,1]$ extends to  a proper
projection $\barr\pr_{h, \mbfx}: \barr \scrm_{h, \mbfx} \to [0,1]$;
\item $\barr\pr_{h, \mbfx}$ is smooth on every stratum of $\barr \scrm_{h,
\mbfx}$.  
\end{itemize}
More precisely, every stratum of $\barr \scrm_{h, \mbfx}$ consists of pairs $(C,
t)$ such that $C$ is possibly {\slsf reducible and not reduced} $h(t)
$-holomorphic curve in the homology class $[\Sigma]$ passing through $\mbfx$. 
Thus every $C$ is a formal sum $C = \sum_i m_i C_i$ of closed irreducible
$h(t)$-holomorphic curves with positive integer multiplicities $m_i$, such that
$[\Sigma] = \sum_i m_i [C_i]$ and $x_1, \ldots x_k \in \supp(C) = \cup_i C_i$. 
The strata are indexed by obvious combinatorial data: homology classes, genera,
multiplicities of single components, and the distribution of the points $x_1,
\ldots x_k$ on the components. The smooth structure on the strata describes
deformation of components in terms of solutions of the equation \eqqref(2.2). 
The topology on the whole compactification $\barr \scrm_{h, \mbfx}$ is the
{\slsf cycle topology} in which every curve $C = \sum_i m_i C_i$ is considered
as a {\slsf closed 2-current} on $X$, see below for details. We refer to
\cite{Sh} for more details on the structure of $\barr \scrm_{h, \mbfx}$.

\newdefi{def2.3a} Let $C_n$ be a sequence of pseudoholomorphic curves in a
manifold $X$ with parameterizations $u_n: S_n \to X$. It converges to a
pseudoholomorphic curve $C^*$ with a parameterization $u^*: S^* \to X$ in the
{\slsf cycle topology} if
\begin{itemize}
\item[\slsf CT1] the boundaries $\d S_n$ and $\d S^*$ have the same number of
circles; moreover, there exists diffeomorphisms $\phi_n : \d S^* \to \d S_n$
such that the maps $u_n \scirc \phi_n: \d S^* \to X$ converge to $u^*\ogran_{\d
S^*}: \d S^* \to X$ in the $C^1$-topology;
\item[\slsf CT2] for any continuous 2-form $\psi$ on $X$ the integrals
$\int_{u_n( S_n)} \psi$ converge to $\int_{u^*( S^*)} \psi$;
\item[\slsf CT3] curves $C_n$ and $C^*$ are holomorphic with respect to almost
complex structures $J_n$ and $J^*$ on $X$ respectively, such that $J_n$ converge
to $J^*$ in the $C^0$-topology.
\end{itemize}
In fact, in the assertions below we shall have even a little bit finer version
of the cycle topology. Namely, the convergence of the structures $J_n \lrar J^*$
will be in the H\"older $C^{0,\alpha} $-topology with some $0< \alpha <1$ except
sufficiently small neighborhoods of the singular points of $C^*$.

\end{defi}

Using the saddle property {\slsf S2)\.} one can show that under condition
$c_1(X)\cdot [\Sigma] >k$ there exists a continuous piecewise smooth section
$s(t) = (C_t, t)$ section of $\barr \pr_{h, \mbfx}: \barr \scrm_{h, \mbfx} \to
[0,1]$. One would obtain the desired symplectic isotopy if one manages to
``push'' such a section $s$ into $\scrm_{h, \mbfx} ^\circ$, \ie deform $s$ into
a section $s'(t)$ with values in $\scrm_{h, \mbfx} ^\circ$, or even in $\scrm_{h,
\mbfx}$. To understand whether such a deformation exists one needs a description
how different strata of $\barr \scrm_{h, \mbfx}$ are attach to each other. Thus
we are led to the question of description of possible symplectic isotopy classes
of nodal curves in a neighborhood of a given singular pseudoholomorphic curve
$C^*$. This question is often related to as the {\slsf local symplectic isotopy
problem}. 

As it was noticed in \cite{Sh}, {\slsf Main Theorem} provides a sufficiently
complete solution of local symplectic isotopy problem for nodal curves in a
neighborhood of a {\slsf reduced} pseudoholomorphic curve $C^*$, \ie in the case
when every irreducible component of $C^*$ is non-multiple. For a precise
statement we need a generalization of some notions for the case of
pseudoholomorphic curves.

\newdefi{def2.4} Let $X$ be a {\sl 4}-manifold, $J_0$ an almost complex structure
on $X$, and $C_0$ a $J_0$-holomorphic curve with a parameterization $u_0: S \to
X$.  Assume that $C_0$ has no multiple component and that the boundary $\d C$ is
non-singular or empty. 

An {\slsf equigeneric deformation} $C_t$ of $C_0$ is given by a family $J_t$ of
almost complex structures on $X$ and a family $u_t: S \to X$ of
$J_t$-holomorphic maps such that $C_t = u_t(S)$ and such that every $u_t$ is an
imbedding near the boundary $\d S$. We assume that the structures $J_t$ and the
parameterization maps $u_t$ depend continuously on $t$. Every pseudoholomorphic
curve $C_1$ which appears in this way is also called an equigeneric deformation
of $C_0$.

A {\slsf maximal nodal deformation} of $C_0$ is a {\slsf nodal} curve $C_1$
which is an equigeneric deformation of $C_0$. As in the usual holomorphic case,
every singular point $p$ of $C_0$ ``splits'' under maximal nodal deformation
into certain number of nodes. This number is called the {\slsf (virtual) nodal
number} of $C_0$ at $p$ and denoted usually by $\delta(C^*, p)$. The sum
$\delta(C^*) \deff \sum \delta(C^*, p_i)$ over all singular points of $C^*$ is
the maximal number of nodes which can be obtained by a deformation of $C^*$
which is small in the cycle topology.

A {\slsf nodal deformation} of $C_0$ is given by a family $J_t$ of almost
complex structures on $X$ and a family $C_t$ of {\slsf nodal} $J_t$-holomorphic 
curves such that 

\begin{itemize}
\item the structures $J_t$ depend continuously on $t$;
\item the curves $C_t$ depend continuously on $t$ with respect to the cycle
topology;
\item every $C_t$ is imbedded near the boundary $\d C_t$; moreover, the
boundaries $\d C_t$ depend continuously on $t$ with respect to the $C^1$-topology.
\end{itemize}

As in the holomorphic case, every small deformation $C_1$ of a {\slsf nodal}
curve $C_0$ is nodal again; however, some nodes of $C_0$ disappear and some
persist. We say that $C_1$ is obtained from a nodal curve $C_0$ by {\slsf
smoothing the nodes $p_1, \ldots, p_l$} of $C_0$ if $C_1$
is a small nodal deformation of $C_0$ and the set of nodes which disappear is 
$\{ p_1, \ldots, p_l \}$.
\end{defi}

The following result about the uniqueness of maximal nodal deformation and
smoothing of a prescribed set of nodes is proved in \cite{Sh}.

\newprop{prop2.3} \sli Let $X$ be a {\sl 4}-manifold and $C^*$ be a
pseudoholomorphic curve whose boundary is either empty or smooth imbedded. Then
two sufficiently small maximal nodal deformations $C_0$ and $C_1$ of $C^*$ can
be connected by an isotopy $C_t$ which is close to $C^*$ in the cycle topology. 

\slii Let $X$ be a {\sl 4}-manifold, $C^*$ be a nodal pseudoholomorphic curve
whose boundary is either empty or smooth imbedded, and $\{p_1, \ldots, p_l\}$ a
prescribed subset of the set of nodes of $C^*$. Then two sufficiently small
deformations $C_0$ and $C_1$ of $C^*$ obtained by smoothing the prescribed nodes
$p_1, \ldots, p_l$ can be connected by an isotopy $C_t$ which is close to $C^*$
in the cycle topology.

In both cases, if $C^*$ is $J^*$-holomorphic and the structure $J^*$ is tamed by
a symplectic form $\omega$, then the isotopy $C_t$ can be chosen
$\omega$-symplectic.
\end{prop}

\newsubsection[sym2]{Existence of symplectic isotopy between nodal surfaces}
The first application of {\slsf Main Theorem} is the positive solution of the
local symplectic isotopy problem for nodal pseudoholomorphic curves without
multiple components.

\newthm{thm2.4} Let $X$ be a {\sl 4}-manifold, $J^*$ an almost complex structure
on $X$, and $C^*$ a $J^*$-holomorphic curve. Assume that $C^*$ has no multiple
component and that the boundary $\d C$ is smooth imbedded or empty.

Let $C$ be some nodal deformation of $C^*$ and $C^\dag$ a maximal nodal
deformation of $C^*$, both sufficiently close to $C^*$ in the cycle topology.
Then there exists an isotopy $C_t$ between $C$ and a small deformation $C^\ddag$
of $C^\dag$ obtained by smoothing an appropriate set of nodes of $C^\dag$.
Moreover, the isotopy $C_t$ can be realized sufficiently close to $C^*$

Moreover, if the structure $J^*$ is tamed by a symplectic form $\omega$ on $X$,
then the isotopy $C_t$ can be made $\omega$-symplectic.
\end{thm}

\proof As it was already indicated, the assertion follows from {\slsf Main
Theorem} and the techniques developed in \cite{Sh}, especially in {\slsf
Subsection 6.2}. Let us outline the modifications needed to adapt the
argumentation used there to our situation.

\smallskip\noindent
{\slsf Special case}. Assume that $X$ is the unit ball in $\cc^2$, the
structure $J^*$ is sufficiently close to the standard structure in $\cc^2$, and
$C^*$ has a single singularity at the origin $0 \in B$. 

\smallskip\noindent
{\slsf Preparatory construction}. Performing an appropriate isotopy, one can reduce
the problem to the situation when $C^*$ is holomorphic. The construction of such
an isotopy used in \cite{Sh} applies here with minor modification.

\smallskip\noindent
{\slsf Induction by complexity of singularities}. In \cite{Mi-Wh}, Micallef and
White has proved that the local behavior of pseudoholomorphic curves is
essentially the same as the one of genuine holomorphic curves. In particular,
one obtains well-defined notions of the topological type of the singularity and
of  codimension of a singularity of a given topological type. The latter is
the codimension of the space of curves with the singularity of the given
topological type in the whole space of curves. 

A parameter version of the result of Micallef and White was proved in \cite{Sh},
{\slsf Section 3}. In particular, the actual codimension of the set of
pseudoholomorphic curves with a singularity of a given topological type is the
expected one, see \cite{Sh} for details. Inductively, we may assume that the
assertion of the theorem holds for all pseudoholomorphic curves whose
singularities have smaller codimension than that of $C^* \subset B$.

\smallskip\noindent
{\slsf Main construction}. One tries to find an isotopy $C_t$ between $C \ddef
C_0$ and a holomorphic curve $C_1$ controlling the behavior of $C_t$ near the
boundary so that $C_t$'s remain close to $C^*$. It is proved in \cite{Sh},
{\slsf Subsection 6.2}, that there exists an isotopy $C_t$ such that
\begin{itemize}
\item $C_t$ is parameterized by $t \in [0,t^+)$ and remains close to $C^*$;
\item for some increasing sequence $t_n$ converging to $t^+$ the sequence
$C_{t_n}$ converges to a curve $C^+$;
\item the curve $C^+$ either is holomorphic in the usual sense or has
singularities of codimension strictly smaller than that of $C^* \subset B$.
\end{itemize}
If the obtained curve is $C^+$ is holomorphic, then the assertion of the theorem
for the special case of a single singularity follows from {\slsf Main Theorem}. 
Otherwise, the assertion  follows by induction.

\smallskip\noindent
{\slsf General case}. One performs appropriate constructions in a neighborhood
of every singular point of $C^*$ and then extend the obtained local families of
deformations of $C^*$ to a global family $C_t$. Such a family $C_t$ can be made
$J_t$-holomorphic since there are no integrability condition on the structures
$J_t$.
\qed

\medskip
Our second application is the positive solution of the (global) symplectic
isotopy problem for nodal surfaces of lower genus. 

\newthm{thm2.5} \sli Let $(X, \omega)$ be a $\cp^2$ with the Fubini-Study
form. Then every two symplectic nodal irreducible surfaces $\Sigma_0, \Sigma_1$
of the same degree and the same genus $g \le 4$ are symplectically isotopic.

\slii Let $X$ be $\cp^2$ blown-up at one point, and $\omega$ a symplectic form
on $X$. Then every two symplectic nodal irreducible surfaces $\Sigma_0,
\Sigma_1$ of the same homology class and the same genus $g \le 2$ are
symplectically isotopic.

\sliii Let $X$ be $S^2 \times S^2$ and $\omega$ a product symplectic form on 
$X$.  Then every symplectic nodal irreducible surface $\Sigma$ of genus $g \le
3$ is symplectically isotopic to an algebraic curve. In particular, there exist
finitely many symplectic isotopy classes of nodal irreducible surface $\Sigma$
of a given genus $g \le 3$ in a given homology class on $S^2 \times S^2$.
\end{thm}

The general idea of the proof is as follows. 
In all three cases there exists the standard complex structure $J\st$ on
$X$ tamed by the symplectic form $\omega$. This means that $(X, J\st)$ is
isomorphic to $\cp^2$, the blown-up $\cp^2$, or $\cp^1 \times \cp^1$,
respectively. We shall show that every symplectic nodal surface $\Sigma \subset
X$ satisfying the hypotheses of the theorem is symplectically isotopic to a
$J\st$-holomorphic curve. The uniqueness of the symplectic isotopy class in the
case of the (blown-up) $\cp^2$ will follow then from the irreducibility of the
Severi variety $V_g(X, [\Sigma])$ of irreducible nodal $J\st$-holomorphic curves
in $X$ of genus $g$ in the homology class $[\Sigma]$. This result is proved by
Harris \cite{Ha} for $\cp^2$ and by Ziv Ran \cite{Ran} for $\cp^2$ blown-up at
one point.

\medskip
Now let $\Sigma \subset X$ be a symplectic nodal surface satisfying the
hypotheses of the theorem. In particular, $\Sigma$ is irreducible and has genus
$g$ at most $4$, $3$, or $2$ according to $X$. To find a symplectic isotopy between
$\Sigma$ and a $J\st$-holomorphic curve we repeat the construction which was
used in \cite{Sh}, {\slsf Subsection 6.3}, and exposed in \refsubsection{sym1}.

\smallskip
First we establish possible values of the ``anti-canonical degree'' $c_1(X)\cdot
[\Sigma]$ for nodal symplectic surfaces satisfying the hypotheses of
\refthm{thm2.5}.

\newlemma{lem2.6} Let $\Sigma$ be a nodal symplectic surface in a symplectic
{\sl4}-fold $(X, \omega)$. Then ``anti-canonical degree'' $c_1(X)\cdot [\Sigma]$
is at least $1$ if $X$ is the blown-up $\cp^2$ and $\Sigma$ is $J$-holomorphic
for some structure $J$ which can be included in a generic 1-parameter family
of structures $J_t$; at least $2$ if $X$ is $S^2 \times S^2$; and at least $3$
if $X$ is $\cp^2$.

Moreover, if the equality holds then, according to the case, $\Sigma$ is
\begin{enumerate}
\item an exceptional sphere,  if $X$ is the blown-up  $\cp^2$;
\item a ``horizontal'' or ``vertical'' line representing the homology class
$[S^2 \times \pt]$ or $[\pt \times S^2]$, respectively, if $X$ is $S^2 \times
S^2$;
\item a ``line'' \ie a sphere of degree $1$, if $X$ is $\cp^2$.
\end{enumerate}
\end{lem}

\state Remark. The same assertion holds in the case when $\Sigma$ is an
algebraic curve in the (blown-up) $\cp^2$ or $\cp^1 \times \cp^1$,
respectively. This classical result follows also from the proof of the lemma.

\proof {\slsf Case $X= \cp^2$}. In this case $(X, \omega)$ is symplectomorphic
to $\cp^2$ equipped with some positive multiple of the Fubini-Study form
$\omega_{FS}$. The group $\sfh_2(\cp^2, \zz)$ is $\zz$ and every $\omega
$-symplectic nodal surface must have positive degree $d$. Then $c_1(X)\cdot
[\Sigma] = 3d$. The genus formula for symplectic nodal surfaces insures that 
$\Sigma$ is an imbedded sphere in the case $d=1$.

\medskip
{\slsf Case $X$ is the\/ $\cp^2$ blown-up at one point}. We use basic properties
of symplectic blown-up in dimension 4 and symplectic exceptional spheres, see
\eg \cite{McD-3}. Assume that $J$ is a generic $\omega $-tame almost complex
structure on $X$ and $\Sigma \subset X$ an irreducible nodal $J$-holomorphic
curve. Furthermore, we assume that $\Sigma$ is not an exceptional sphere since
otherwise $c_1(X)\cdot [\Sigma] =1$. Then that there exists a $J$-holomorphic
exceptional sphere $E \subset X$. It follows from the genericity of $J$ that $E$
meets $\Sigma$ only at smooth points and transversally. Perturbing $J$, we can
make $J$ integrable near $E$. Denote $d_E \deff [\Sigma] \cdot [E]$. Then $d_E$
is a non-negative integer. Contracting $E$ we obtain
\begin{itemize}
\item a compact 4-manifold $X'$ diffeomorphic to $\cp^2$;
\item a point $p_E$ which appears instead of the exceptional sphere $E$, such
that $X'\bs \{ p_E \}$ is canonically identified with $X \bs E$;
\item the symplectic form $\omega'$ on $X'$ whose restriction on $X'\bs \{ p_E
\}$ coincides with $\omega\ogran_{X \bs E}$;
\item an $\omega'$-tame almost complex structure $J'$ on $X'$ which is
integrable near $p_E$;
\item a $J'$-holomorphic curve $\Sigma'$ in $X'$ such that $\Sigma'\bs \{ p_E
\}$ coincides with $\Sigma \bs E$ and such that $\Sigma'$ has $d_E$ non-singular
transversal branches at $p_E$.
\end{itemize}
Note that $\Sigma'$ is irreducible since $\Sigma$ is assumed to be so. Denote
by $d$ the degree of $\Sigma'$ in $\cp^2$.  Then the homology class of $\Sigma$
is $[\Sigma] = d L - d_E E$ where $L$ denotes a ``line in $X$'', \ie the lift to
$X$ of a generic $J'$-holomorphic line in $X'$. In particular, $c_1(X)\cdot
[\Sigma] = 3d -d_E$.

We assert that $d_E \le d$ and the equality holds \iff $d_E = d =1$. Indeed,
perturbing $\Sigma'$ at $p_E$ we obtain a nodal symplectic surface with
$\frac{d_E(d_E -1)} {2}$ new nodes instead the singularity of $\Sigma'$ at
$p_E$. The genus formula for this perturbation reads
\[
g(\Sigma) = g(\Sigma') = \frac{(d-1)(d -2)}{2} - \delta(\Sigma) - 
\frac{d_E(d_E -1)}{2}  
\]
where $\delta(\Sigma)$ is the number of nodes of $\Sigma$. This implies the
desired inequality $d_E \le d$ and shows that the equality holds in the unique
case $d_E = d =1$. This case corresponds to the $J'$-holomorphic line in $X'$
passing through $p_E$.

Now, the inequality $d_E \le d$ together with the formula $c_1(X)\cdot [\Sigma]
= 3d -d_E$ yield the desired inequality $c_1(X)\cdot [\Sigma]\ge 2$. 

\state Remark. Observe that as the consequence of the argumentation above we
obtain that the equality $c_1(X)\cdot [\Sigma] = 2$ holds in the unique case
when $\Sigma$ is an imbedded sphere with trivial normal bundle meeting the
exceptional sphere $E$ at a single point. This means that $\Sigma$ is a fiber of
a {\slsf $J$-holomorphic ruling on $X$}, see \cite{McD-2} and \cite{McD-Sa-2}
for details.

\medskip
{\slsf Case $X= S^2 \times S^2$}. In this case $\sfh_2(X, \zz) \cong \zz \oplus
\zz$. We use the ``almost complex'' geometry of ruled symplectic 4-manifold, see
\cite{McD-2} and \cite{McD-Sa-2} for details. It provides the existence of an
$\omega$-tame almost complex structure $J$ on $X$ with the following properties:
\begin{itemize}
\item $\Sigma \subset X$ is a $J$-holomorphic curve;
\item there exist $J$-holomorphic curves $L_h$ and $L_v$ which represent the
``horizontal'' and ``vertical'' homology classes $[S^2 \times \pt]$ and $[\pt
\times S^2]$, respectively.
\end{itemize}
It follows that $[\Sigma] = a [L_h] + b [L_v]$ with {\slsf non-negative}
integers $a= [\Sigma] \cdot [L_v]$ and $b= [\Sigma] \cdot [L_h]$, and that
$c_1(X) \cdot [\Sigma] = 2a + 2b$. Thus $c_1(X) \cdot [\Sigma] \ge 2$, and the
equality holds \iff $\Sigma$ is either ``horizontal'' or ``vertical'' line as
above. \qed

\bigskip
Turn back to the proof of \refthm{thm2.5}. Recall that $J\st$ denotes an
$\omega$-tame integrable structure such that $(X, J\st)$ is isomorphic to
$\cp^2$ or $\cp^1 \times \cp^1$ or the blown-up $\cp^2$ according to the case we
have. Find an $\omega$-tame almost complex structure $J_0$ making $\Sigma$ a
$J_0 $-holomorphic curve, denoted by $C_0$. Set $k \deff c_1(X) \cdot [\Sigma]
-1$. Fix $k$ distinct points $\mbfx=(x_1, \ldots, x_k)$ on $C_0$. Perturbing
$C_0$ and the points, we may assume that $x_1,\ldots, x_k$ are in general
position \wrt the structure $J\st$ in the following sense. For any closed
surface $S$, not necessary connected, the moduli space $\scrm_{J\st, \mbfx}(S,
X, [\Sigma])$ of $J\st$-holomorphic (and hence {\sl algebraic}) curves of the
homology class $[\Sigma]$ with normalization $S$ passing through $\mbfx$ is
either empty or a complex space of the expected dimension.

Fix a generic path $h(t)$ of $\omega$-tame almost complex structures $J_t 
\deff h(t)$ connecting $J_0$ with $J\st=J_1$. Without loss of generality we may 
assume that $J_t$ depend $C^\ell$-smoothly on $x\in X$ and $t$ for some
$\ell\gg0$. Our hope is to find an isotopy $C_t$ between $\Sigma= C_0$ and a
$J_1$-holomorphic curve which consists of $J_t$-holomorphic curves. Trying to
construct such a family $C_t$ for maximal possible interval we obtain

\newprop{prop6.3.1} There exists a $t^+ \in (0,1]$ which is {\sl maximal} \wrt 
the following condition: 

For any $t <t^+$ there exists a $J_t$-holomorphic curve $C_t$ such that 
\begin{itemize}
\item[\slip] $C_t$ passes through the fixed points $\mbfx=(x_1,\ldots, x_k)$;
\item[\sliip] $C_t$ is non-multiple and irreducible;
\item[\sliiip] the curve $C_0$ is symplectically isotopic to the curve 
obtained from some maximal nodal deformation $C'_t$ of $C_t$ by smoothing an
appropriate set of nodes of $C'_t$.
\end{itemize}
\end{prop}

\smallskip
Let $t_n$ be an increasing sequence converging to $t^+$. Fix $J_{t_n}
$-holomorphic curves $C_n$ with these properties. Property \sliii implies that
the $C_n$ have the same homology class as $C_0$. Going to a subsequence we may
assume that they converge to a $J_{t^+}$-holomorphic curve $C^+$ in the cycle
topology.

\newprop{prop6.3.2} Under the hypotheses of \refthm{thm2.5}, assume that $C^+$
has multiple components. Then $C^+$ has two irreducible components, $C'$ of
multiplicity $1$ and $L$ of multiplicity $2$ such that, according to the case,
\begin{enumerate}
\item $C'$ has genus $2$ and $L$ is an exceptional line, if 
$X$ is the blown-up $\cp^2$;
\item $C'$ has genus $3$ and $L$ is a horizontal or vertical line, if 
$X$ is $S^2 \times S^2$;
\item $C'$ has genus $4$ and $L$ is a line, if $X$ is $\cp^2$.
\end{enumerate}
Moreover, the curve $L \cup C'$ is nodal and the marked points $\mbfx$ are
disjoint from the nodes of $L \cup C'$.
\end{prop}

The latter condition means that $L \cup C'$ is in generic position and there are
no further degeneration or incidences than those stipulated by the hypotheses of
the proposition.

\medskip
\proof Let $C^+ = \sum m_i C^+_i$ be the decomposition of $C^+$ into irreducible
components with multiplicities $m_i$. Set $\mu \deff c_1(X) \cdot [C^+] = c_1(X)
\cdot [\Sigma]$ and $\mu_i \deff c_1(X) \cdot [C^+_i]$. Let $g_i$ be the
(geometric) genus of $C_i$ and $k_i$ the number of the marked points $\mbfx$
lying on $C_i$. It follows then that $k_i \le \mu_i + g_i -1$. The reason is
that otherwise the expected dimension of the space of irreducible curves of
genus $g_i$ in the homology class $[C_i]$ passing through $k_i$ points is
negative; hence the existence of such a constellation would contradict the
condition of the generality of $h(t)$, see \eg {\slsf Subsection 2.4} of
\cite{Sh}. Besides, we have the obvious (in)equalities $\mu = \sum \mu_i$, $\sum
k_i \ge k= \mu -1$ and $g\deff g(\Sigma) \le \sum m_i g_i$. Taking into account
the inequality $\mu_i \ge 3$, $2$, or $1$, according to the cases of
\lemma{lem2.6}, and distinguishing the case of equality, we see that multiple
components are possible only in the cases described in the proposition.

The genericity properties of $C'$ follows from the condition of the genericity of
$\mbfx$ and $h(t)$. Namely, similarly to the usual holomorphic (and hence
algebraic) case, every additional incidence or degeneration condition, such as
appearance of a cusp or a triple point, makes the expected dimension of the
corresponding constellation negative, which would again contradict the
genericity, see \cite{Sh}. \qed

\medskip
Let us distinguish the cases according to the structure of the curve $C^+$.

\medskip\noindent
{\slsf Case 1. $C^+$ is irreducible.} We claim that $t^+=1$ in this case.
Assuming the contrary it is sufficient to show that for some $t^{++}>t^+$ there
exists a $J_{t^{++}}$-holomorphic curve $C^{++}$ with the properties given in
\propo{prop6.3.1}. To do this we fix some parameterization $u^+: S^+ \to C^+
\subset X$ and consider the relative moduli space $\scrm_{h, \mbfx}(S^+, X)$ of
$J_t= h(t) $-holomorphic curves which are parameterized by $S^+$, pass through
$\mbfx$ and lie in the homology class $[\Sigma]$. This space is non-empty
because it contains $C^+$. It follows the from the results of \cite{Sh},
especially {\slsf Subsection 4.5}, that for some $t^{++} >t^+$ such a curve
$C^{++}$ does exist. 

Now, since $t^+=1$, the structure $J_{t^+}$ is $J\st$, the standard one, and
$C^+$ is an irreducible algebraic curve in $(X, J\st)$. Let $g^+$ be the
geometric genus of $C^+$. It follows now from {\slsf Proposition 2.1} of
\cite{Ha} that every component of the variety $V(|C^+|, g^+)$ of irreducible
curves of geometric genus $g$ in $|C^+|$ is of expected dimension $c_1(X)\cdot
[C^+] + g^+ -1$ and contains a nodal curve. Consequently, $C^+$ can be included
in a 1-dimensional family $\{ C_\lambda \}$ whose generic member $C_\lambda$ is
an irreducible nodal curve of the same genus as $C^+$. Observe that such
$C_\lambda$ is a maximal nodal deformation of $C^+$. Then, smoothing an
appropriate set of nodes of $C_\lambda$, we obtain the desired algebraic curve
which is symplectically isotopic to $C_0 = \Sigma$ by \refthm{thm2.4}. 
This yields the proof of \refthm{thm2.5} for the special {\slsf Case
1} of irreducible $C^+$. The existence of the desired smoothing is provided by

\newlemma{lem2.9} Let $X$ be a non-singular complex projective surface and $C
\subset X$ a nodal curve without multiple components such that $c_1(X)\cdot C_i$
is positive for every irreducible component $C_i$ of $C$. Then every prescribed
set of nodes of $C$ can be smoothed by some deformation of $C$.
\end{lem}

\proof Let $\wt C$ be the normalization of $C$, $u: \wt C \to X$ the induced
immersion, $\scri_C$ the defining ideal of $C \subset X$ and $\scrn_C \deff
(\scri_C / \scri_C ^2) ^*$ the normal sheaf of $C$. Then there exists a natural
projection map $p: \scro_X(TX) \to \scrn_C$ with the following properties:
\begin{itemize}
\item the kernel $\ker(p)$ is naturally isomorphic to the sheaf $\scro_C(TC)$ of
sections of the tangent bundle of (the normalization of) $C$;
\item the image $\im(p)$ is naturally isomorphic to the sheaf
$\scro_C(TX/du_*(TC))$ of sections of the normal bundle $N_C \deff TX/du_*(TC)$;
\item the cokernel $\scrn_C /\im(p)$ is isomorphic to the sum $\sum_i
\scro_{x_i}$ over all nodal points $x_i$ of $C$.
\end{itemize}
More precisely, we construct appropriate sheaves on the normalization $\wt C$
and then push them forward onto $C$ or $X$ by means of $u$.

Now let $\mbfx=\{x_1, \ldots, x_k\}$ be some set of nodes of $C$. Denote by
$\scrn_{C, \mbfx}$ the sheaf on $C$ which coincides with $\scrn_C$ at each
smooth point of $C$ and each of the nodes $\{x_1, \ldots, x_k\}$, and with the
image $\im(p)= \scro_C(N_C)$ at each remaining node. The deformation theory (see
\eg \cite{Pal-1} and \cite{Pal-2}) insures that
\begin{itemize}
\item the space of deformations of $C$ which smooth the prescribed nodes $\mbfx$
is given by a Kuranishi model $\Phi: B \to \sfh^1(C, \scrn_{C, \mbfx})$ for some
holomorphic map $\Phi$ defined in some ball $B$ in $\sfh^0(C, \scrn_{C, \mbfx})$;
\item the natural projection $\sfh^0(C, \scrn_{C, \mbfx}) \to \sfh^0(C, 
\sum_{i=1}^k \scro_{x_i})$ describes the smoothing of nodes. In particular, a
deformation with the tangent vector $v\in \sfh^0(C, \scrn_{C, \mbfx})$ smoothes
the node $x_i$ \iff the projection of $v$ in $\sfh^0(C, \scro_{x_i})$ does not
vanish.
\end{itemize}

Let $C_i$ be an irreducible component of $C$, $g_i$ its geometric genus, and
$N_{C_i}\deff TX/ du_*(TC_i)$ the corresponding normal bundle. Then $c_1(
N_{C_i}) = c_1(X) \cdot [C_i] + (2g_i-2) > 2g_i-2$ by the hypothesis of the
lemma.  Consequently, $\sfh^1(C_i, \scro(N_{C_i})) =0$ for each single normal
bundle. Thus the obstruction group $\sfh^1(C, \scrn_{C, \mbfx})$ vanishes and
every prescribed set of nodes $\mbfx$ can be smoothed. \qed

\medskip\noindent
{\slsf Case 2. $C^+$ is reducible but without multiple components.} 
Let $C^+_i$ be the irreducible components
of $C^+$. Then the (bi)degree of each $C^+_i$ is strictly less than the 
(bi)degree of $C^+$. Applying induction, we may assume that the assertion of 
\refthm{thm2.5} holds for every $C^+_i$. Moreover, we may also suppose that
for $t\in [t^+,1]$ there exist families $\{C^+_{i,t}\}$ of $J_t$-holomorphic 
curves with the following properties:
\begin{itemize}
\item $C^+_{i,t^+} = C^+_i$, \ie every family $\{C^+_{i,t}\}$ starts from
$C^+_i$ at $t^+$;
\item for every $t\in [t^+,1]$ the curve $C^+_t \deff \cup_i C^+_{i,t}$
is nodal;
\item smoothing appropriate set of nodes on $C^+_t = \cup_i C^+_{i,t}$ we obtain
a curve which is symplectically isotopic to $C_0 = \Sigma$.
\end{itemize}

For the final value $t=1$, the existence of {\slsf nodal} curve $C^+_1= \cup_i
C^+_{i,1}$ with the desired properties follows from {\slsf Proposition 2.1} of
\cite{Ha}.

In particular, smoothing appropriate set of nodes on the ``final''
curve $C^+_1= \cup_i C^+_{i,1}$ gives the desired algebraic curve which is
symplectically isotopic to $C_0=\Sigma$.

\medskip
It remains to consider 

\medskip\noindent
{\slsf Case 3. $C^+$ has multiple components.} Recall that $C^+$ was obtained as
the limit of a sequence of $J_{t_n}$-holomorphic curves $C_{t_n}$. To simplify
notation, we write $J^+$ instead of $J_{t^+}$, $C_n$ instead of $C_{t_n}$, and 
$J_n$ instead of $J_{t_n}$.

Notice that the limit $C_n \lrar C^+$ is understood in the cycle
topology. However, we obtain more information about the behavior of $C_n$
near $C^+$ if we take the limit in the {\slsf stable map topology} instead of
the cycle one.

For the definition of the stable map topology and related notions in full
generality we refer to {\slsf Section 5} of \cite{Sh} and \cite{Iv-Sh}, as also
to \cite{Ha-Mo} and \cite{Fu-Pa} for the algebraic setting. In our setting, the
limit object is given by an abstract closed nodal curve $\wh C^+$ equipped with
$J^+ $-holomorphic map $u^+: \wh C^+ \to X$ which have the following properties:
\begin{enumerate}
\item[\slsf St1)] Let $\{ \wh C_\bfla \}$ be any semi-universal family of
deformations of $\wh C^+$ such that $\wh C_{\bfla^+}$ is the curve points $\wh
C^+$ itself. Then there exists a sequence of parameters $\bfla_n$ converging to
$\bfla^+$ such that, after going to a subsequence, $\wh C_{\bfla_n}$ is
isomorphic to the normalization $\wt C_n$ of $C_n$.
\item[\slsf St2)] the image $u^+(\wh C^+)$, counted with multiplicities, is
$C^+$. 
\item[\slsf St3)] If $\wh C_i^+$ is a rational irreducible component of $\wh
C^+$ and the number of nodal points on $\wh C_i^+$, counted with multiplicities,
is less than $3$, then $u^+$ is non-constant on $\wh C_i^+$.
\end{enumerate}
The first condition means that there we can imbed $\wh C^+$ and the
normalizations $\wt C_n$ in $X \times \cp^N$, pseudoholomorphic \wrt the
structures $J^+ \times J_{\cp^N}$ and $J_n \times J_{\cp^N}$ respectively, so
that the images $\wt C_n \subset X \times \cp^N$ will converge to $\wh C^+
\subset X \times \cp^N$ in the cycle topology, and so that the projection of
these images onto $X$ gives the sequence $C_n$ converging to $C^+$. In
particular, the map $u^+$ can be obtained as the projection from $\wh C^+
\subset X \times \cp^N$ onto $C^+ \subset X$, and the second condition follows.
The last condition excludes the appearance of redundant components and insures
the uniqueness of the limit in the stable map topology.  Observe also that by
the second condition the {\slsf arithmetic} genus of $\wh C^+$ is the {\slsf
geometric} genus of $C_n$.

The crucial point in treating of {\slsf Case 3} is study of the deformation
problem of the pair $(\wh C^+, u^+)$ in the stable map topology. We start with
establishing the possibilities for the structure of $(\wh C^+, u^+)$. Obviously,
we must have a component $\wh C'$ mapped by $u^+$ onto the component $C'$ of
$C^+$ as in \propo{prop6.3.2}.  Denote by $\wh C''$ the remaining part of $\wh
C^+$.

\newlemma{lem2.10} Under the hypotheses of \refthm{thm2.5} and \propo{prop6.3.2}, 

\noindent
\sli $\wh C'$ is the normalization of $C'$ and $u^+: \wh C' \to C'$ is the
normalization map;

\noindent
\slii there are the following possibilities for the remaining part $\wh C''$:

\begin{enumerate}
\item[(A)] $\wh C''$ consists of two rational components $\wh C''_1$ and $\wh C''_2$,
each mapped by $u^+$ isomorphically onto the line $L$ and attached to $\wh C'$
at points $z_1^\times, z_2^\times \in \wh C'$, respectively; the images
$u^+(z_1^\times)$ and $u^+(z_2^\times)$ are two distinct intersection points of
$C'$ and $L$.

\smallskip
\item[(B1)] $\wh C''$ is rational and attached to $\wh C'$ at a point
$z^\times_1$ whose image $u^+(z_1^\times)$ is an intersection point of $C'$
and $L$; the map $u^+ : \wh C'' \to X$ is a too shitted covering of $L \subset
X$ branched over two distinct points $y_1, y_2 \in L$.

\smallskip
\item[(B2)]  $\wh C''$ consists of two rational components $\wh C''_1$ and
$\wh C''_2$, each mapped by $u^+$ isomorphically onto the line $L$; $\wh
C''_1$ is attached to $\wh C'$ at a point $z_1^\times\in \wh C'$ and $\wh
C''_2$ to $\wh C''_1$ at a point $z_2^\times \in \wh C''_1$; the image
$u^+(z_1^\times)$ is an intersection point of $C'$ and $L$; the image
$u^+(z_2^\times)$ lies on $L$ apart from $u^+(z_1^\times)$.

\smallskip
\item[(B3)] $\wh C''$ consists of three rational components $\wh C''_0$, $\wh
C''_1$, and $\wh C''_2$; $\wh C''_0$ is attached to $\wh C'$ at a point
$z^\times_0\in \wh C'$; $\wh C''_1$ and $\wh C''_2$ are attached to $\wh
C''_0$ at two distinct points $z^\times_1, z_2^\times \in \wh C''_0$, which
are distinct also from $z^\times_0$; $u^+$ maps $\wh C''_1$ and $\wh C''_2$
isomorphically onto $L$ and $\wh C''_0$ constantly into the point
$u^+(z^\times_0)$ which is an intersection point of $C'$ and $L$.

\end{enumerate}
\sliii If $X$ is the blown up $\cp^2$ and $L$ is an exceptional line, then
only case (A) is possible.
\end{lem}

\proof The first assertion follows by comparing the geometric genera of $C'$
and $\wh C'$. The same argument implies that the remaining part $\wh C''$ must
consist of trees of rational curves. Thus $L$ can be covered either by one or by
two distinct rational curves. Elementary combinatorics shows that the cases
(A) and (B1--B3) are the only possibilities for such trees of rational curves.

Now assume that $X$ is the blown up $\cp^2$ and the sequence $C_n$ converges to
one of the the constellations (B1--B3). Then we can choose an appropriate piece
$C^\circ _n$ of each $C_n$ such that $C^\circ _n$ are connected and the limit of
$C^\circ _n$ in the cycle topology consists of the exceptional line $L$ with
multiplicity $2$ and a disc $D$ transversal to $L$. The intersection index of
$C_n ^\circ$ with $L$ must be $[C_n ^\circ] \cdot [L] = [D] \cdot [L] + 2 \cdot
[L]^2 =-1$. Now observe that $C_n ^\circ$ are holomorphic with respect to
structures $J_n$ converging to $J^+$ such that there exists a sequence of $J_n
$-holomorphic exceptional lines $L_n$ which converges to the line $L$. 
Consequently, $[C_n ^\circ] \cdot [L] = [C_n ^\circ] \cdot [L_n] \ge 0$. With
this contradiction the proof is finished.
\qed

\state Remark. Observe that the constellations (A) and (B3) are rigid \ie
determined by the curve $C^+$ and the combinatorics. To the contrary, we obtain
moduli in the constellations (B1) and (B2), namely, positions of the branching
points $y_1$ and $y_2$ in the second constellation, and position of the point
$u^+(z_2 ^\times)$ in the third one. The constellation (B1) degenerates in (B2)
as $y_1$ and $y_2$ collapse apart from $u^+(z_1^\times)$, and in (B3) as $y_1$
and $y_2$ collapse with $u^+(z_1^\times)$. These combinatorial data and varying
parameters is the additional information we obtain taking the limit in the
stable map topology instead of the cycle one.

\medskip
Trying to deform $(\wh C^+, u^+)$ in the stable map topology into an irreducible
curve we come to the {\slsf gluing problem} for pseudoholomorphic curves. Let us
resume the results of {\slsf Subsection 5.3} of \cite{Sh} on this topic which we
shall use.

\newdefi{def2.6} A {\slsf pants} $P$ is a complex curve which can be obtained from
$\cp^1$ by removing 3 disjoint discs with smooth boundary. {\slsf Boundary
annuli} in a pants $P$ are disjoint annuli $A_1, A_2, A_3 \subset P$ each
adjacent to some boundary circle of $P$

The {\slsf standard smoothing of a node} is the family 
\[
\scra_\lambda \deff \{ (z_1, z_2) \in \Delta^2: z_1 \cdot z_2 =\lambda \}
\]
with the parameter $\lambda$ varying in a disc $\Delta(\eps) \deff \{ |\lambda|
< \eps \}$ of radius $\eps <1$. It deforms the {\slsf standard node} $\scra_0$,
consisting of two discs $\Delta_1$ and $\Delta_2$ with the canonical coordinates
$z_1$ and $z_2$ respectively, into annuli $\scra _\lambda$, $\lambda \ne0$. The
{\slsf boundary annuli} $A_1, A_2 \subset \scra _\lambda$ are given by
\[
A_1 \deff \{ (z_1, z_2) \in \scra_\lambda : 1-\delta < |z_1| <1\}
\qquad
A_2 \deff \{ (z_1, z_2) \in \scra_\lambda : 1-\delta < |z_2| <1\}
\]
with $\delta < \frac{1-\eps}2$.  We consider $A_1$ and $A_2$ with the canonical
coordinates $z_1$ and $z_2$, respectively, as a ``constant'' part inside
deforming curves $\scra_\lambda$.

For an almost complex manifold $(X, J)$, we denote by $\scrp(\scra_\lambda, X,
J)$ the space of $J$-holomorphic maps $u: \scra_\lambda \to X$ which are
$C^1$-smooth up to boundary. In the case of $\scra_0$ such a map $u: \scra_0 \to
X$ is given by its components $u_1: \Delta_1 \to X$ and $u_2: \Delta_2 \to X$,
both $J$-holomorphic, such that $u_1(0) = u_2(0)$. For any compact (nodal) curve
$C$ with the smooth boundary $\d C$, possibly empty, the space $\scrp(C, X, J)$
is defined in a similar way.
\end{defi}

\newprop{prop2.11} \sli For any compact (nodal) curve $C$ without closed
components the space $\scrp(C, X, J)$ has a natural structure of a Banach
manifold. 

\slii For any given structure $J^*$, a compact (nodal) curve $C$ without closed
components, and a map $u^* \in \scrp(C, X, J^*)$ there exists an open
neighborhood $\scru \subset \scrp(C, X, J^*)$ of $u^*$ and a map $G=G(u, J):
\scru \to \scrp(C, X, J)$ depending smoothly on $u\in \scru$ and on a structure
$J$ sufficiently $C^1$-close to $J^*$ such that, for $J$ fixed, the
map $G_J: \scru \to \scrp(C, X, J)$ is an open smooth imbedding.

\sliii The restriction maps $R_\lambda: \scrp(\scra_\lambda, X, J) \to \scrp(A_1,
X, J) \times \scrp(A_2, X, J)$ given by $R_\lambda(u) \deff (u\ogran_{A_1},
u\ogran _{A_2})$ are smooth closed imbeddings. 

\sliv For any given $J^*$ and $u^* \in \scrp(\scra_0, X, J^*)$ there exists an
open neighborhood $\scru \subset \scrp(\scra_0, X, J^*)$ of $u^*$ and a map
$G=G(\lambda, u, J): \scru \to \scrp(\scra_\lambda, X, J)$ defined for $u\in
\scru$, $\lambda$ sufficiently close to $0$, and for structures $J$ sufficiently
$C^1$-close to $J^*$, such that:
\begin{itemize}
\item $G$ is continuous in $\lambda$ and $C^1$-smooth in $u$ and $J$;
\item for $\lambda$ and $J$ fixed, the map $G: \scru \to \scrp(\scra_\lambda, X,
J)$ is an open $C^1$-smooth imbedding.
\end{itemize}
Moreover, in cases \sliii and \sliv the $C^1$-smoothness is uniform in
$(\lambda, u, J)$.
\end{prop}

The meaning of the last part of the proposition is that we can ``glue'' the
components $u^*_{1,2}$ of any given pseudoholomorphic map $u^*: \scra_0$ into a
pseudoholomorphic map $u: \scra_\lambda \to X$, also varying the almost complex
structure. To apply the proposition in our situation we must decompose $\wh C^+$
into appropriate pieces. For the proof of the following assertion we refer to
\cite{Iv-Sh}.

\newprop{prop2.12} There exist a covering $\{V_a \}$ of $\wh C^+$ and families
of deformations $V_{a, \lambda_a}$ of some pieces $V_a$ with the following
properties:
\begin{itemize}
\item Every piece $V_a$ is isomorphic to the standard node $\scra_0$, or the
disc $\Delta$, or an annulus $\scra_\lambda$, or a pants.
\item Each intersection $V_a \cap V_b$, if non-empty, is an annulus $A_{ab}$
which is a boundary annulus for both $V_a$ and $V_b$.
\item The pieces included in the deformation families are all nodal pieces $V_a
\cong \scra_{\lambda^+_a=0}$ and some annular pieces $V_a \cong \scra_{
\lambda^+_a \ne 0}$. The deformation family for such a piece $V_a \cong
\scra_{\lambda^+_a}$ is of the form $V_{a, \lambda_a} = \scra_{\lambda_a}$ with
$\lambda_a$ varying in a small neighborhood of $\lambda^+_a$.
\item Let $\bfla$ be the system of all $\lambda_a$'s which appear as the
parameter of the families $V_{a, \lambda_a}$'s, and let $\wh C_\bfla$ be the
curve obtained by replacing each varying piece $V_a$ by the piece $V_{a,
\lambda_a}$. Then $\{ \wh C_\bfla \}$ is a semi-universal family of deformations
of $\wh C^+$.
\end{itemize}
\end{prop}

We divide the obtained parameters $\bfla=(\lambda_1, \ldots, \lambda_l)$ into
two groups: $\bfla''=(\lambda''_1, \ldots, \lambda''_{l''})$ each describing the
smoothing of the corresponding node on $\wh C^+$, and the remaining $\bfla'=
(\lambda'_1, \ldots, \lambda'_{l'})$, $l' + l'' =l$. Thus we obtain
$\bfla''=(\lambda''_1, \lambda''_2)$ in the cases (1) and (3) of \lemma{lem2.10},
$\bfla''=(\lambda''_1)$ in the case (2), and $\bfla''=(\lambda''_0, \lambda''_1,
\lambda''_2)$ in the case (4). Let $\bfla^+= (\bfla'{}^+, \bfla''{}^+ )$ be the
set of parameters corresponding to the curve $\wh C^+$ so that $\bfla''{}^+ =0$.

\smallskip
Using the covering $\{ V_a \}$ we describe the problem of deformation of $(\wh
C^+, u^+)$ in terms of compatibility of deformations of single pieces $V_a$ and
the restrictions of $u^+$ onto $V_a$'s. Namely, let us fix a small $C^1
$-neighborhood $\scru_J$ of $J^+$, small neighborhoods $\scru_{\bfla'}$ and
$\scru_{\bfla''}$ of $\bfla'{}^+$ and $\bfla''{}^+$ in the spaces of parameters 
$\bfla'$ and $\bfla''$ respectively, and, for each $V_a$, a small
neighborhood $\scru_a$ of the restriction $u^+_a \deff u^+ \ogran_{V_a}$ in the
space $\scrp(V_a, X, J^+)$. Consider the map
\[
\textstyle
\scrg: \left(\prod_a \scru_a \right) 
\times  \scru_{\bfla'} \times  \scru_{\bfla''} 
\times  \scru_J \lrar
\prod_{a\ne b} \scrp(A_{ab}, X, J),
\]
where the product $\prod_{a\ne b} \scrp(A_{ab}, X, J)$ is taken over all pairs
$(a,b)$ for which the intersection $V_a \cap V_b$ is a non-empty annulus
$A_{ab}$. For such a pair $(a,b)$, the component $\scrg_{ab}$ of $\scrg$ is
defined as follows. We take the $a$-th component $u_a$ of $\mbfu \in \prod_a
\scru_a$, compute its deformation $u'_a \deff G(u_a, \lambda_a, J)$ or $u'_a
\deff G(u_a, J)$ according to the type of $V_a$, the obtained map $u'_a$ lies in
$\scrp(V_{a, \lambda_a} J)$ or $\scrp(V_a, J)$ respectively, and then restrict
$u'_a$ onto $A_{ab}$.

Observe that every annulus $A_{ab}$ appears twice, as $V_a \cap V_b$ and as
$V_b \cap V_a$, but the components $\scrg_{ab}$ and $\scrg _{ba}$ do not
coincide in general. Moreover, the set of conditions
\[
\scrg_{ab}(\mbfu, \bfla, J) = \scrg _{ba}(\mbfu, \bfla, J)
\qquad
\text{for each pair $(a,b)$}
\]
is the compatibility condition on the pieces $G(u_a, \lambda_a, J)$ or $G(u_a,
J)$ to be the restrictions on $V_a$ of a well-define a $J$-holomorphic map $u:
\wh C_\bfla \to X$. Thus, denoting by $\scrd_J \subset \prod_{a\ne b}
\scrp(A_{ab}, X, J)$ the ``diagonal set'' given by the set of conditions $u_{ab}
= u_{ba}$, we obtain the set-theoretic equality 
\[
\scrp(\wh C_\bfla, X, J) = \scrg(\cdot, \bfla, J)\inv(\scrd_J),
\]
which holds locally near $(\wh C^+, J^+)$. Observe also that $\scrg$ is only
continuous in $\bfla''$ but still $C^1$-smooth in the remaining variables
$\mbfu$, $\bfla'$, and $J$. 

\newlemma{lem2.13} Let $\bfla^*$, $J^*$, and $u^* \in \scrp(\wh C_{\bfla^*}, X,
J^*)$ be close to $\bfla^+$, $J^+$, and $u^+$ respectively. Set $u^*_a \deff u^*
\ogran_{V_{a, \lambda_a^*}}$ and $\mbfu^* \deff (u^*_a) \in \prod_a \scru_a$.

Then the map $\scrg( \mbfu, \bfla', \bfla''{}^*, J^*)$, with the arguments
$\mbfu$ and $\bfla'$ varying and $\bfla''{}^*$ and $J^*$ fixed, is transversal
to the sub\-mani\-fold $\scrd_{J^*}$ at the point $(\mbfu^*, \bfla^*, J^*)$. 
\end{lem}

\proof The transversality means that the image of differential of the map
$\scrg( \cdot, \cdot, \bfla''{}^*, J^*)$ at the point $(\mbfu^*, \bfla^*, J^*)$
is the whole normal space to $\scrd_{J^*} \subset \prod_{a\ne b} \scrp(A_{ab},
X, J^*)$ at $\scrg( \mbfu^*, \bfla^*, J^*)$. An equivalent assertion is that the
deformation problem described by $\scrg( \cdot, \cdot, \bfla''{}^*, J^*)$ is
unobstructed because the cokernel of the differential in question in the the
normal space to $\scrd_{J^*}$ is the obstruction space to the deformation
problem.

We may assume that $\bfla^*= \bfla^+$, $\wh C^* =\wh C^+$, $J^* =J^+$, and $u^*
= u^+$. The general case follows from this special one by the following
argument. A surjective linear {\slsf Fredholm} map between Banach spaces remains
surjective after a small perturbation. We compute the deformation problem in two
steps as follows: first, we consider the deformation problems for each component
$\wh C'$ and $\wh C''_i$ of $\wh C^+$, and then impose the conditions of
``attaching''.

\smallskip
{\slsf Step 1}. Observe that the parameters $\bfla'$ parameterize a complete
family of deformations of $C'$. This follows from the fact that $\wh C^+$
differs from $C'$ by trees of rational curves. Consequently, the map $\scrg(
\cdot, \cdot, \bfla''{}^+, J^+)$ describes the problem of deformation of $C'$ as
a {\slsf parameterized} $J^+$-holomorphic curve of the given geometric genus $g
= g(C')$. Observe also that the curve $C'$ is immersed and $c_1(X)\cdot [C']
>0$. These two conditions imply that the deformation problem is unobstructed,
see \eg \cite{H-L-S} or \cite{Sh}, {\slsf Section 2}. The same argument applies
for the components $\wh C''_i$.

\smallskip
{\slsf Step 2}. After solving the problems of the first step, we obtain local
deformations families of $J^+$-holomorphic maps: $u'_{s'}: S \to X$, defined on
a closed real surface $S$ of genus $g=g(C')$, and $u''_{i, s_i}: S^2 \to X$, one
for each component $\wh C'$ and $\wh C''_i$, respectively. To fit together in a
map of a {\slsf connected} curve $\wh C_{\bfla', \bfla'' {}^+}$, the maps
$u'_{s'}$ and $u''_{i, s_i}$ must satisfy certain ``attaching conditions''
defined as follows. Each nodal point $z^\times_i$ on $\wh C_{\bfla', \bfla''
{}^+}$ has two pre-images on the components $\wh C'_{\bfla'}$ and $\wh C''_i$,
say $z^+_i$ and $z^-_i$, and the images of these points in $X$ must coincide. 
The transversality of this ``attaching problem'' is equivalent to the original
transversality. For this purpose possibility to move arbitrarily the image of
one of the points $z^+_i$ and $z^-_i$ is sufficient. The latter condition is
equivalent to the transversality of the problems of deformations of the curves
$C'$ and $\wh C''_i$ constrained by the condition of passing through given
points.

We contend that this new deformation problem is unobstructed. Let us consider
the special case when the curve $\wh C^+$ is as in the case (A) of \lemma{lem2.10}
and $L$ is an exceptional line. In this case the component $C'$ of $C^+$ meets
$L$ at two points at least. This implies that $c_1(X) \cdot [C'] \ge3$ since
otherwise $C'$ would meet $L$ at a single point, see the remark in the proof of  
\lemma{lem2.6}. Now, since $C'$ is immersed and $c_1(X) \cdot [C']$ is strictly
larger than the number $k=2$ of the constraining points, the problem of
deformation of $C'$ constrained at $k=2$ points is unobstructed. This yields the
desired transversality for the special case we consider. 

The other cases can be treated similarly.
\qed

\medskip
As a corollary of \lemma{lem2.13} we obtain the local symplectic isotopy in a
neighborhood of $C^+$. 

\newcorol{cor2.14} \sli Let $(\wh C^+, u^+)$ be as in the case {\sl (A)} of
\lemma{lem2.10} and $(\wh C_0, u_0)$, $(\wh C_1, u_1)$ two small deformations of
$(\wh C^+, u^+)$ in the stable map topology, such that $C_i \deff u_i(\wh C_i)$
are irreducible and nodal. Then there exists a symplectic isotopy $C_t$ between
$C_0$ and $C_1$ close to $C^+$ in the cycle topology.

\slii Let $(\wh C^+_0, u^+_0)$ and $(\wh C^+_1, u^+_1)$ be as in the cases {\sl
(B1--3)} of \lemma{lem2.10} and $(\wh C_0, u_0)$, $(\wh C_1, u_1)$ two small
deformations of $(\wh C^+_i, u^+_i)$ in the stable topology, $i=0,1$
respectively. Assume that $C_i \deff u_i(\wh C_i)$ are irreducible and nodal. 
Then there exists a symplectic isotopy $C_t$ between $C_0$ and $C_1$ close to
$C^+$ in the cycle topology.
\end{corol}

 Observe that the almost complex structure can also vary.

\proof Let $J^*$ be a structure close to $J^+$. Set $\scrm_{J^*}
\deff \cup_{\bfla''} \scrg(\cdot, \cdot, \bfla'', J^*) \inv (\scrd _{J^*})$ and
let $\scrm\sing_{J^*}$ be the set of parameters $(u, \bfla, J^*) \in
\scrm_{J^*}$ where $u(\wh C_\bfla)$ is not nodal and irreducible. It follows
from \lemma{lem2.13} that $\scrm_{J^*}$ is a {\slsf topological} manifold in a
neighborhood of $(\mbfu^+, \bfla^+, J^+)$ and that $\scrm\sing_{J^*}$ has
Hausdorff codimension $\ge2$ in $\scrm_{J^*}$.
This fact and \lemma{lem2.13} imply part \sli of the corollary.

For part \slii we use an additional possibility to connect $(\wh C^+_0,
u^+_0)$ and $(\wh C^+_1, u^+_1)$ by a path $(\wh C^+_t, u^+_t)$ continuous in
the stable map topology such that $u^+_t(\wh C^+_t)$ is constantly $C^+$.
\qed

\medskip
Now we are ready to finish

\statep Proof. of \refthm{thm2.5}. Recall that it remains to consider the
following situation: There exists a sequence $C_n$ of pseudoholomorphic nodal
curves such that each $C_n$ is symplectically isotopic to $\Sigma$ and such that
there exists the limit $(\wh C^+, u^+)$ of $C_n$ in the stable map topology. 
Furthermore, the possibilities for the structure of $(\wh C^+, u^+)$ are given
by \lemma{lem2.10}. \lemma{lem2.13} and \refcorol{cor2.14} insure the
possibility of restoration of the symplectic isotopy class of $\Sigma$ by $C^+ =
u^+(\wh C^+)$ and the combinatorial data.

The scheme of the proof is the same as before: First, we show that there exists
a symplectic isotopy $C^+_t$ between $C^+_0 \deff C^+$ and a holomorphic curve
$C^+_1$, and then deform $C_1^+$ into a holomorphic nodal curve in the
symplectic isotopy class of $\Sigma$. 

Proving the existence of the desired symplectic isotopy $C^+_t$ we apply the
induction in the ``anti-canonical degree''. Namely, by \lemma{lem2.6} we have
$c_1(X) \cdot [C'] < c_1(X) \cdot [\Sigma]$ for the component $C'$. Thus there
exists a symplectic isotopy $C'_t$ between $C' = C'_0$ and a holomorphic curve
$C'_1$. The existence of a similar symplectic isotopy for $L$ is follows
directly from the following fact: For a generic path of tame structures $J_t$
and a generic choice of points $x_1, \ldots, x_k$ with $k\deff c_1(X) \cdot [L]
-1$ there exists a {\slsf unique} path $L_t$ formed by $J_t$-holomorphic curves
in the homology class $L$. This fact was exploited by several authors, see \eg
\cite{Bar}. It follows then that both isotopies $C'_t$ and $L_t$ can be made
$J_t$-holomorphic for the same path of tamed structures $J_t$. Then for a
generic choice of isotopies $C'_t$ and $L_t$ the curves $C^+_t \deff C'_t \cup
L_t$ will form the desired symplectic isotopy.

The combinatorial data are translated along the path $C^+_t$ onto the obtained
holomorphic curve $C_1^+$. Since \lemma{lem2.13} and \refcorol{cor2.14} hold
also for the structure $J\st$, we can deform $C^+_1$ into a nodal $J\st
$-holomorphic curve in the symplectic isotopy class of $\Sigma$. \refthm{thm2.5}
follows. \qed

\bigskip

\ifx\undefined\bysame
\newcommand{\bysame}{\leavevmode\hbox to3em{\hrulefill}\,}
\fi

\def\entry#1#2#3#4\par{\bibitem[#1]{#1}
{\textsc{#2 }}{\sl{#3} }#4\par\vskip2pt}


\end{document}